\documentclass[twocolumn]{autart}    

\usepackage{cite}

\usepackage{amsmath}
\usepackage{amsmath,amssymb,amsfonts} 
\usepackage{graphicx}

\usepackage{algorithm}
\usepackage{algpseudocode}

\usepackage{textcomp}
\usepackage{epstopdf}
\usepackage{booktabs} 
\usepackage{amsfonts}       
\usepackage{multicol}
\usepackage{multirow} 
\usepackage{subcaption}
\usepackage{xcolor} 
\graphicspath{ {images/} } 
\newcommand{\norm}[1]{\lVert{#1}\rVert}

\newcommand{\bmat}[1]{\begin{bmatrix}#1\end{bmatrix}}

\newcommand{\mcl}[1]{\mathcal{#1}}
\newcommand{\mbf}[1]{\mathbf{#1}} 

\newcommand{\N}{\mathbb{N}}
\newcommand{\R}{\mathbb{R}}

\definecolor{UBCgreen}{rgb}{0.04706, 0.55, 0.03} 
\definecolor{UBCgreen2}{rgb}{0.04706, 0.40, 0.13} 

\begin{document}

\begin{frontmatter}

\title{Picard Iteration for Parameter Estimation in Nonlinear Ordinary Differential Equations {  using Low-Quality Data}} 

\thanks[footnoteinfo]{This work was supported in part by the National Science Foundation under Grants NSF CCF-2323532 and NSF/NIH 2054354 awarded
 to M. M. Peet. (Corresponding author: Aleksandr Talitckii.)}

\author[ASU]{Aleksandr Talitckii}\ead{atalitck@asu.edu},    
\author[ASU]{Matthew M. Peet}\ead{mpeet@asu.edu},               

\address[ASU]{School for the Engineering of Matter, Transport and Energy, Arizona State University, Tempe, AZ, 85298 USA}  


\begin{abstract}                          
We consider the problem of using experimental time-series data for parameter estimation in nonlinear ordinary differential equations, {  focusing on the case where the data is noisy, sparse, irregularly sampled, includes multiple experiments, and does not directly measure the system state or its time derivative.} To account for such low-quality data, we propose a new framework for gradient-based parameter estimation which uses the Picard operator to reformulate the problem as constrained optimization with infinite-dimensional variables and constraints.
{  
We then formulate the Karush-Kuhn-Tucker (KKT) conditions necessary for optimality and define a convergent sequence of approximations to these KKT conditions obtained by replacing the solution map by the $n$-th order Picard iterate. Then, for any element of this sequence, and by exploiting the contractive properties of the Picard operator, we propose a gradient-contractive algorithm which (under regularity and convexity assumptions) is guaranteed to converge to a solution of these approximated KKT conditions.
} 
Finally, the algorithms are then tested on a battery of models and a variety of datasets in order to demonstrate robustness and improvement over alternative approaches. 

\end{abstract}

\end{frontmatter}

\section{Introduction}
Mechanical, electrical, chemical and biological processes are often dynamic with highly nonlinear response to disturbance, actuation, and initial state. Regulation of such processes requires accurate models of this nonlinear behaviour. Ideally, the dynamics of the process can be modelled through the combination of well-established physical principles. In some cases, however (e.g. biological systems), the physical principles governing the behaviour are speculative, over-simplified, or only partially understood. Furthermore, even when the physics of the problem are well-understood, there may still be substantial uncertainty in parameters of the model due to natural variation, errors in measurement, oversimplification of the physics or difficulty in direct measurement of these parameters. In such cases, there may be only a few uncertain parameters in the model, and they may be constrained to lie in some set (e.g. positivity is a common form of parameter constraint).

When it is not possible to measure the parameters of the system directly, these parameters must be inferred indirectly by examination of the response of the system to variations in input and initial state~\cite{ljung1998system}. In the extreme case, where there is no understanding of the physics, we have black-box system identification~\cite{brunton2022data} --- e.g. Koopman operators~\cite{yeung2019learning, mauroy2016linear}; extended dynamic mode decomposition (EDMD)~\cite{williams2015data}; sparse identification of non-linear dynamics (SINDy)~\cite{kaheman2020sindy, naozuka2022sindy, brunton2016discovering}; autoregressive models~\cite{billings2013nonlinear}; and neural networks~\cite{raissi2019physics}.

When the physics are only partially understood or the system parameters are constrained, however, we have grey-box system identification~\cite{bohlin2006practical} -- a more challenging problem than black-box modeling.
Specifically, suppose the system to be identified is known to have the form\vspace{-3.5mm}
\begin{align}
\dot x(t)&= f(t, x(t), \theta) \qquad x(0) = x_0\notag \\
y(t) &= g(t, x(t), \theta) \label{eqn:nonlinearODE} \\[-8mm]\notag
\end{align}
where the functions $f$ and $g$ are given, $x(t) \in \R^{n_x}$ is the internal state, $x_0$ is the initial state, and $y(t)\in \R^{n_y}$ is some measurable output. The vector, $\theta\in \Theta \subset \R^{n_\theta}$, represents the unknown parameters in the model. The dependence of $f$ and $g$ on time, $t$, is typically used to represent the effect of known inputs, $u(t)$ --- in which case we might equivalently write the functions as $f(x(t), u(t), \theta)$ and $g(x(t), u(t), \theta)$. 

 { Of course, when measurements of $\dot x(t)$ are available, or $x(t)$ is measurable ($g=f$ or $g(t,x)=x$, respectively) and the sampling instances are sufficiently dense, it is relatively straightforward to use least squares or some other form of regression to estimate the parameters, $\theta$~\cite{ljung1998system}. However, when sampling instances are sparse in time or only output measurements are available, the system identification problem becomes significantly more challenging~\cite{brunton2022data}.   

Specifically, in this paper, we consider the case of grey-box parameter estimation based on noisy irregular sampling of system output, $y(t)$, and where the data may be generated using multiple experiments, each with different initial condition $x_{0,i}$ and sampling instances. Note that implicit in this framework is that the initial condition, $x_{0,i}$ for each experiment is unknown and must be estimated along with the system parameters, $\theta$ (similar to the inverse problem~\cite{aster2005parameter}).}
 %
%

For a given parameterized vector field, $f(t,x,\theta)$, we suppose that Eqn.~\eqref{eqn:nonlinearODE} is well-posed and define the corresponding solution map, $\phi_{f}(t, x, \theta)$, as the unique function which satisfies
\vspace{-3.5mm}
\[
\partial_t \phi_{f}(t, x, \theta) = f(t, \phi_{f}(t, x, \theta), \theta),\quad 
\phi_{f}(0, x, \theta) = x \vspace{-3.5mm}
\] for all $t \ge 0$, $x \in \R^{n_x}$ and $\theta\in \Theta$. While it is almost never possible to obtain an analytic expression for the solution map of a nonlinear system, the $\phi_{f}$ notation allows us to efficiently represent the problem of data-based system identification. Specifically, suppose that experimental data is available from solutions of Eqn.~\eqref{eqn:nonlinearODE} consisting of measurements $y_i$ at a sequence of discrete times $\{t_i\in\R\}_{i=1}^{N_s}$ so that $y_i = (1+n_i) g(t_i, \phi_{f}(t_i, x, \theta), \theta) + m_i$ 
where $\theta,x$ are unknown and $m_i\in \R^{n_y}$ and $n_i \in \R$ represent measurement errors.  We denote the set of all such measurements as $Y = \{(y_i, t_i) \hspace{-0.5mm}\in \hspace{-0.5mm}\R^{n_y}\hspace{-0.5mm}\times\hspace{-0.5mm}\R \; : \; i = 1, ..., N_s\}.$
%

We may now formulate the parameter estimation/system identification problem as \vspace{-3.5mm}
\begin{equation}
\min_{\substack{\theta\in\Theta, x\in\R^{n_x}}}  L_{ls}(x, \theta, \phi_{f}(\cdot, x, \theta)), \label{eqn:optimization}\vspace{-3.5mm}
\end{equation}
where for simplicity we use the least-squares \textit{Loss Function}~\cite{bard1974nonlinear} defined as \vspace{-6mm}
\begin{equation*}
 L_{ls}(x, \theta, \phi_{f}(\cdot, x, \theta))\hspace{-0.5mm}=\hspace{-0.5mm}\frac{1}{2N_s} \sum_{i=1}^{N_s} \|y_i - g(t_i, \phi_{f}(t_i, x, \theta), \theta)\|_2^2. \vspace{-3.5mm}
\end{equation*}

%
{ 
Problem~\eqref{eqn:optimization} is an unconstrained optimization problem~\cite{ljung2010perspectives} in parameters $x,\theta$. However, for nonlinear systems, the inclusion of the nonlinear solution map, $\phi_f$, typically makes the problem non-convex -- implying that descent-based algorithms will obtain at best local optimality. Furthermore, because explicit expressions for the solution map $\phi_f$ are rarely available, obtaining even locally optimal solutions to this problem is challenging.  
}

One approach to solving Problem~\eqref{eqn:optimization} is to use numerical simulation in place of the solution map~\cite{biegler1986nonlinear} and then apply gradient-free, black box optimization tools such as numerical gradient-approximation~\cite{byrd1988approximate}, parameter cascade~\cite{ramsay2007parameter, hooker2016collocinfer}, nonparametric estimators~\cite{brunel2008parameter, chen2008efficient, cao2012penalized}, 
 genetic algorithms~\cite{duong2002system, sivanandam2008genetic}, and simulated annealing~\cite{bertsimas1993simulated}. Unfortunately, however, such methods are computationally inefficient and may fail to converge to even locally optimal solutions.

By contrast, gradient-based methods are guaranteed to converge to local optima but require some method for
computing the gradient of the solution map, $\nabla_{x, \theta} \phi_{f}$. 
Because analytic expressions for the solution map are rarely available, however, such a gradient descent approach is not practical.  
 
Faced with the limitations of black-box optimization and lack of analytic solution maps, and motivated by the use of gradient descent to obtain local optima, we consider now an alternative formulation of the parameter estimation/system identification problem posed in Eqn.~\eqref{eqn:optimization}, which does not require knowledge of the solution map. Specifically, in Sec.~\ref{sec:problem_statement}, we reformulate the unconstrained optimization problem as the equivalent constrained optimization problem defined as \vspace{-3.5mm}
\begin{align}   \label{eqn:intro_optimization_problem_constrained}
   \min_{\substack{x \in X, \theta \in \Theta, \mbf u \in C(\R)}} \hspace{-4mm} &\; L_{ls}(x, \theta, \mbf u)\\[-3mm]
     \text{s.t.} &\;  {\mbf u}(t) = x+\int_{0}^t f(s,\mbf u(s),\theta)ds \qquad \forall t \ge 0. \notag\\[-8mm] \notag
\end{align}
In Optimization Problem~\eqref{eqn:intro_optimization_problem_constrained}, we have replaced the solution map $\phi_f$ by the variable $\mbf u$ and added the constraint that $\mbf u$ be a solution of the given system. Clearly, Problems~\eqref{eqn:optimization} and~\eqref{eqn:intro_optimization_problem_constrained} are equivalent.
Moreover, the constraint, ${\mbf u}(t) = x+\int_{0}^t f(s,\mbf u(s),\theta)ds$, may be conveniently represented using the Picard operator as $\mbf u=\mcl P_{x, \theta} \mbf u $, where \vspace{-5mm}
\begin{equation*}
	(\mcl P_{ x, \theta} \mbf u)(t):=  x + \int_{0}^t f(s, \mbf u(s),  \theta) ds.\vspace{-3.5mm}  
\end{equation*} 
 Although Optimization Problem~\eqref{eqn:intro_optimization_problem_constrained} does not require the solution map, it introduces an infinite-dimensional variable $\mbf u$ and associated equality constraint which must hold for all times, $t$, on which the data is defined. Solving such problems without an explicit representation of the feasible set is known to be hard~\cite{polyak1987introduction}. 
%

To address these challenges, in Sec.~\ref{sec:algorithms} we propose a new class of gradient-contractive algorithms (Alg.~\ref{alg3:modified_picard}) for solving the formulation of the parameter estimation problem in Eqn.~\eqref{eqn:intro_optimization_problem_constrained}. In this algorithm, each iteration has two steps -- gradient computation and contraction to the feasible set. Specifically, the first step computes the gradient of the loss function and updates the initial condition and parameter variables as $x_{k+1},\theta_{k+1}$ where \vspace{-3.5mm}
\begin{align*}
 x_{k+1} &= x_{k} - \alpha \nabla_{x} L_{ls}(x_k, \theta_k, \mcl P_{ x_k, \theta_k}^n \mbf u_{k})\\
 \theta_{k+1} &= \theta_k - \alpha \nabla_{\theta} L_{ls}(x_k, \theta_k, \mcl P_{ x_k, \theta_k}^n \mbf u_{k}), \\[-8mm]
\end{align*}
and where the $n^{th}$ order Picard iteration, $\mcl P_{ x, \theta}^n$, is introduced into the objective so that the gradient accounts for multipliers  corresponding to the KKT conditions for the constraint $\mbf u=\mcl P_{x, \theta} \mbf u $ and where the step size $\alpha > 0$ is computed through a line search. The second step contracts the variable $\mbf u_{k}(t)$ to the new feasible set (determined by $x_{k+1}$ and $\theta_{k+1}$) using the contractive property of the Picard iteration as\vspace{-3.5mm}
\[
	\mbf u_{k+1}(t) = (1- \sigma) \mbf u_{k}(t) + \sigma (\mcl P_{ x_{k+1}, \theta_{k+1}} \mbf u_k)(t) \vspace{-3.5mm}
\] 
for some $\sigma \in (0, 1]$. This second step translates the solution closer to the infinite-dimensional manifold described by the feasible set of Problem~\eqref{eqn:intro_optimization_problem_constrained} -- imitating the projection step of gradient-projection algorithms and resolving the problems introduced by the infinite-dimensional nature of the equality constraints. This approach also eliminates the need to explicitly parameterize the infinite-dimensional variable, $\mbf u$. { A limitation of this approach, of course, is the need for the Picard operator to be contractive on the interval of time on which the data was sampled. However, as described in Sec.~\ref{sec:algorithms}, an extended form of Picard operator may be used, which is contractive over arbitrary intervals (assuming well-posedness of the solution map), but requires the inclusion of additional decision variables for each such extension.

To establish convergence of the proposed class of algorithms to local optimality, we may impose regularity and convexity on the objective function and vector field. However, for nonlinear vector fields, the convexity condition will likely fail for $n\ge 1$. Hence, we necessarily limit our analysis to local convergence to locally optimal solutions. Such a limitation is not particularly strong, however, since many parameter identification problems in nonlinear systems may not have a unique solution or such parameters may be locally but not globally identifiable~\cite{barreiro2023origins}.}

Note that Picard iteration has previously been used for analysis of black-box system identification~\cite{kunze1999solving}, and parameter estimation in~\cite{haeseler2018parameter,slavov2020picard}. Specifically, in~\cite{haeseler2018parameter} gradient descent of the unconstrained Problem~\eqref{eqn:optimization} was considered for signal transduction networks, wherein a Picard iterate was used in lieu of the solution map (See Alg.~\ref{alg1:gradient_algorithm}). Unfortunately, convergence of the Picard iteration is typically limited to short time intervals, thus limiting the generality of this approach. An attempt to extend this approach to longer time intervals using multiple shooting methods was proposed in~\cite{slavov2020picard}. However, convergence of neither method has been proven or extended to sparse data, irregular sampling times or measured outputs.


In contrast to the work in~\cite{haeseler2018parameter,slavov2020picard}, the use of a constrained formulation of the parameter estimation problem, the use of gradient contractive algorithms, and the proposed extended form of Picard iteration (Alg.~\ref{alg4:ext_picard}) used in this paper allow for establishment of convergence proofs, arbitrary sampling intervals, and measured outputs. To establish these results, in Sec.~\ref{sec:main_algorithm}, we first consider the general class of two-step gradient contractive algorithms and establish conditions for convergence based on properties of the contractive operator, $\mcl P$. Next, in Sec.~\ref{sec:picard}, we examine the Picard operator $\mcl P_{ x, \theta}$ and establish conditions for convergence to the solution map and its gradient. In Sec.~\ref{sec:conditions_for_Algs}, we show these results apply to the extended Picard operator (Alg.~\ref{alg4:ext_picard}) and provide conditions for convergence of Alg.~\ref{alg4:ext_picard} to a local solution of the parameter estimation problem. The results are then applied to a rigorous battery of numerical tests, including: Van der Pol oscillator to evaluate the effect of the regularization parameter; the FitzHugh-Nagumo Neuron to evaluate the effect of irregular sampling; the Rosenzweig-MacArthur predator-prey model to evaluate the effect of sparse data; a tumor growth model to investigate the question of identifiability; and the Lorenz system for comparison with SINDy and black-box system optimization methods.

\subsection*{Notation}
We use $\N$, $\R$, and $\R_{+}$, to denote the natural numbers, the real numbers, and the non-negative real numbers, respectively.  For a given $J \in \N$ we denote $\overline{1, J} = \{j \in \N \;|\; 1 \leq j\leq J\}$. For $a, b \in \R_+$ we denote ${\lfloor a/b \rfloor} :=\underset{m\in\N, a \geq mb}{\texttt{argmax}} m$ and $a \bmod{b} := a- \lfloor a/b\rfloor b$.
For a given sets $X, Y$ we denote $\mcl F(X, Y)$ to be the set of all functions $f:X\rightarrow Y$. For compact $S \subset \R^n$, we denote $C^k(S)$ to be the space of $k$ times continuously differentiable functions, $f \in \mcl F(S,\R^m)$ with norm $\|f\| = \sup_{s\in S} \|f(s)\|_2$. We use $C_{a}(S)\subset C(S) := C^0(S)$ to denote the ball of radius $a$ as $C_a(S) = \{\mbf u\in C(S) \; |\; \|\mbf u\| \leq a\}$ and use $C_a$ when the domain is clear from context. 
For $u \in C^1(S)$, we use $\nabla_x u$ to denote the gradient as $\nabla_x u  = \bmat{\frac{\partial}{\partial x_1}u, \dots, \frac{\partial}{\partial x_m}u }^T$.  
For compact $X \subseteq \R^{n}$ we denote $\Pi_X:\R^{n} \rightarrow X$ to be the projection operator so that $\Pi_X(y) = \underset{x \in X}{\texttt{argmin}} \|x- y\|_2$. We say $f\in \mcl F(X,\R)$ is strongly convex on $X\subset \R^n$ with modulus $\mu > 0$ if $f(x) - \frac{\mu}{2} \|x\|^2$ is convex on $X$. We say $f\in \mcl F(X\times Y,\R)$ is strongly convex on $X\subset \R^n$ uniformly in $y$ if there exists $\mu\geq 0$ such that $f(x, y) - \frac{\mu}{2}\|x\|^2$ is convex on $X$ for all fixed $y \in Y\subset \R^m$. For $X\subset \R^n$ and $Y\subset \R^m$, we say that $f\in \mcl F(X , Y)$ is Lipschitz continuous with constant $K$ if for all $x_1, x_2 \in X$ we have $\|f(x_1) - f(x_2)\| \leq K \|x_1 - x_2\|$ and for $Z \subset \R^p$, $f\in \mcl F(X\times Z,Y)$ is Lipschitz continuous with respect to $x$ with constant $K$ if for all $x_1, x_2 \in X$ and $z \in  Z$ we have $\|f(x_1, z) - f(x_2, z)\| \leq K \|x_1 - x_2\|$.

\section{Formulating the Parameter Estimation Problem for Nonlinear Differential Equations}\label{sec:problem_statement}

Consider a parameterized nonlinear differential equation model in state-space representation.   \vspace{-3.5mm}
\begin{align}
\dot {\mbf u}(t) &= f(t, {\mbf u}(t), \theta) \qquad {\mbf u}(0)  = x \notag\\
y(t) & = g(t, {\mbf u}(t), \theta) \label{eqn:ODE} \\[-7mm] \notag
\end{align}
where $t\in \Gamma$ is time (on interval $\Gamma$), ${\mbf u}(t) \in \R^{n_x}$ is system state, $\theta \in \Theta$ are unknown parameters (restricted to compact convex $\Theta \subset \R^{n_\theta}$), $x$ are unknown initial conditions (restricted to compact convex $X \subset \R^{n_x}$), and $y(t) \in  \R^{n_y}$ are measured outputs. The function $f\in C^1(\Gamma \times\R^{n_x}\times\Theta)$ is the vector field and $g\in C^1(\Gamma\times\R^{n_x} \times\Theta)$ is the state-to-output map. We assume existence and uniqueness, and continuous dependence of solutions on initial conditions and parameters, implying existence of a solution map, $\phi:\Gamma\times X\times \Theta \rightarrow \R^{n_x}$, which satisfies \vspace{-3.5mm}
    \begin{equation}
        \partial_t \phi(t, x, \theta) = f(t, \phi(t, x, \theta), \theta), \quad
        \phi( 0, x, \theta ) = x. \label{eqn:solutionmap} \vspace{-3.5mm}
    \end{equation}

Now suppose we are given a set of measurements of the dynamical system for some initial condition $x$ and true parameter value, $\theta$ as \vspace{-3.5mm}
\begin{equation*}
    Y\hspace{-1mm} =\hspace{-1mm}\{ (y_i, \hspace{-0.5mm}t_i) \;\hspace{-1mm} | \hspace{-1mm}\; y_i \hspace{-0.5mm}=\hspace{-0.5mm} (1+n_i)g(t_i, \phi(t_i, x, \theta), \theta) + m_i,\; \hspace{-1mm}
 i \in \overline{1, N_s}  \},\vspace{-3.5mm}
\end{equation*}
where $m_i\in \R^{n_y}$ and $n_i \in \R$ represent measurement errors.

The parameter estimation problem is to find initial condition, $x \in X$, and parameter values, $\theta\in \Theta$, that minimize the least squares error in the predicted output relative to the set of measurements. This can be formulated as either an unconstrained or a constrained optimization problem depending on whether we know the solution map. If the solution map, $\phi$, is known, the problem can be formulated as an unconstrained optimization problem of the form\vspace{-3.5mm}
\begin{equation}
    \min_{\substack{x \in X,  \theta \in \Theta}} L_{ls}(x, \theta, \phi(\cdot, x, \theta)),\vspace{-3.5mm} \label{eqn:optimization_problem_unconstrained}
\end{equation}
where \vspace{-5mm}
\begin{equation*}
L_{ls}(x, \theta, \phi(\cdot, x, \theta)) = \frac{1}{2N_s}\sum_{i=1}^{N_s}\hspace{-1mm} \|y_i - g(t_i, \phi(t_i, x, \theta), \theta)\|^2.\vspace{-3.5mm}
\end{equation*} 
%
Since the solution map is unknown, however, we may reformulate the parameter estimation problem by instead introducing a variable, $\mbf u$, representing the solution of the ODE and add the constraint that this variable be a solution of the ODE. \vspace{-3.5mm} 
\begin{align} \label{eqn:optimization_problem_constrained}
   \min_{\substack{x \in X, \theta \in \Theta, \mbf u \in C(\Gamma)}} \hspace{-4mm} &\; L_{ls}(x, \theta, \mbf u)\\[-3mm] 
     \text{s.t.} &\;  {\mbf u}(t) = x+\int_{0}^t f(s,\mbf u(s),\theta)ds, \notag  \\[-8mm]\notag
\end{align}
where the integral constraint is a necessary and sufficient condition for $\mbf u(t)$ to satisfy the ODE (Eqn.~\eqref{eqn:ODE}) with initial condition $\mbf u(0)=x$ and parameter values $\theta$. The advantage of this formulation is that it does not require knowledge of the solution map. However, the disadvantage is that we have introduced an infinite-dimensional optimization variable in the form of $\mbf u$.

\subsection{Reformulation using the Picard Operator}\label{subsec:reformulation}
The unconstrained optimization problem in Eqn.~\eqref{eqn:optimization_problem_unconstrained} is expressed in terms of the solution map (which is unknown) and the constrained optimization problem in Eqn.~\eqref{eqn:optimization_problem_constrained} is expressed in terms of the variable $\mbf u$ (which is an infinite-dimensional variable and where the constraint must hold at an infinite number of times). The problems with both these formulations may be ameliorated through use of what is known as the Picard operator,  $\mcl P_{x, \theta}$.
{ \begin{defn}\label{def:picard_operator}
For a given $t_0 \in \Gamma$, $T > 0$ and continuous vector field, $f:\Gamma \times\R^{n_x} \times \Theta \rightarrow \R^{n_x}$, we define the Picard operator $\mcl P_{t_0, x, \theta} : C[0,T]\rightarrow  C[0,T] $ \vspace{-3.5mm}
 \begin{equation}
 (\mcl P_{t_0, x, \theta} \mbf u)(t) := x+ \int_{0}^{t} f(s + t_0, \mbf u(s), \theta) ds.\label{eqn:picard_operator_fixed} \vspace{-3.5mm}
 \end{equation} 
 and when $t_0\hspace{-1mm}=\hspace{-1mm}0$, we simplify the notation as $\mcl P_{x, \theta}\hspace{-1mm}:=\hspace{-1mm}\mcl P_{0, x, \theta}$.
 \end{defn}}
 The critical property of the Picard operator is that if $T$ is sufficiently small, iterations of the Picard operator converge to the solution of $\dot x=f(x)$. That is, for {any} $\mbf u(t)$, $\mcl P_{x, \theta} \circ \cdots \circ \mcl P_{x, \theta} \mbf u \rightarrow \mbf u^*$ where $\mcl P_{x, \theta} \mbf u^*=\mbf u^*$.

The Picard operator allows us to \textit{approximate} Optimization Problem~\eqref{eqn:optimization_problem_unconstrained} for fixed $n \in \N$ as \vspace{-3.5mm}
\begin{equation}\label{eqn:optimization_problem_unconstrained_approximated}
    \min_{\substack{x \in X, \theta \in \Theta}} L_{ls}(x, \theta, \mcl P_{ x, \theta}^n \mbf u),  \vspace{-3.5mm}
\end{equation}
 and to \textit{restate} Optimization Problem~\eqref{eqn:optimization_problem_constrained} as\vspace{-3mm}
\begin{equation} \label{eqn:optimization_problem_constrained_approximated}
   \min_{\substack{x \in X, \theta \in \Theta, \\ \mbf u \in C(\Gamma)}} \hspace{-2mm}  L_{ls}(x, \theta, \mbf u)\quad  
     \text{s.t.} \quad    \mcl P_{ x, \theta} \mbf u = \mbf u.\vspace{-3.5mm}
\end{equation} 
Now, the approximation of unconstrained Optimization Problem~\eqref{eqn:optimization_problem_unconstrained} in Problem~\eqref{eqn:optimization_problem_unconstrained_approximated} has eliminated the need for knowledge of the solution map, but introduced a very complicated dependence on initial state, $x$ and parameter, $\theta$. Meanwhile, the restatement of constrained Optimization Problem~\eqref{eqn:optimization_problem_constrained} in~\eqref{eqn:optimization_problem_constrained_approximated} does not seem to have any immediate benefits. As will be discussed in the following section, however, the contractive properties of the solution map will allow us to eliminate the need for infinite-dimensional variables and constraints in~\eqref{eqn:optimization_problem_constrained_approximated}.



\section{Gradient Descent for Parameter Estimation}\label{sec:algorithms}
%

%
\begin{figure}[!t] 
\begin{minipage}[t]{0.5\textwidth}
\begin{algorithm}[H]
\begin{algorithmic}
{\small
\item 
\texttt{INPUT:} $N$ -- the number of iterations \\
\hspace{11mm} $n$ -- the order of Picard iterations \\
\hspace{5mm} $x_0 \in X$, $\theta_0 \in \Theta$, $\mbf u_0 = x_0$, $\alpha  > 0$, { $\varepsilon > 0$, $k = 0$}\\

%
\texttt{REPEAT}\\
\qquad $ \theta_{k+1} := \Pi_\Theta[ \theta_{k} - {\alpha \nabla_\theta L_{ls}(x_{k}, \theta_{k}, \mcl P^n_{x_k, \theta_k} \mbf u_0 )}]$ \\
\qquad $ x_{k+1} :=  \Pi_X[x_{k} - {\alpha \nabla_{x} L_{ls}(x_{k}, \theta_{k}, \mcl P^n_{x_k, \theta_k} \mbf u_0 )}]$ \\ 
\qquad $k = k + 1$,\\
{ \texttt{UNTIL} $\|\nabla_{x, \theta} L_{ls}(x_{k}, \theta_{k}, \mcl P^n_{x_k, \theta_k} \mbf u_0 )\| \leq \varepsilon$ or $k > N$.} \\
\texttt{OUTPUT:} $x_{N}, \theta_{N}$
}
\end{algorithmic}
\caption{Gradient Descent Algorithm}  \label{alg1:gradient_algorithm}
\end{algorithm}
\end{minipage}\vspace{-1mm}

\begin{minipage}[t]{0.5\textwidth}
\begin{algorithm}[H]
\begin{algorithmic}
{\small
\item
 \texttt{INPUT:} $N$ -- the number of iterations \\
 \hspace{5mm} $x_0 \in X$, $\theta_0 \in \Theta$, ${{\mbf u}}_0 = x_0$, $\alpha  > 0, \sigma \in (0, 1]$, { $\varepsilon > 0$, $k = 0$}\\
%
\texttt{REPEAT}\\
%
\qquad $\theta_{k+1} := \Pi_\Theta[ \theta_{k} - {\alpha \nabla_\theta L_{ls}(x_{k}, \theta_{k}, {{\mbf u}}_{k})}]$ \\
\qquad $x_{k+1} :=  \Pi_X[x_{k} - {\alpha \nabla_{x} L_{ls}(x_{k}, \theta_{k}, {{\mbf u}}_{k})}]$ \\
 \qquad${{\mbf u}}_{k+1}(t)  := (1-\sigma) {{\mbf u}}_{k}(t) + \sigma(\mcl P_{x_{k+1}, \theta_{k+1}} {{\mbf u}}_{k})(t)$ \\
 \qquad $k = k + 1$,\\
{ \texttt{UNTIL} $\|\nabla_{x, \theta}  L_{ls}(x_{k}, \theta_{k}, {{\mbf u}}_{k})\| \leq \varepsilon$ or $k > N$.} \\
\texttt{OUTPUT:} $x_{N}, \theta_{N}, {{\mbf u}}_{N}$
}
\end{algorithmic}
\caption{Gradient-Contractive Algorithm}  \label{alg2:picard}
\end{algorithm}
\end{minipage}\vspace{-1mm}

\begin{minipage}[t]{0.5\textwidth}
\begin{algorithm}[H]
\begin{algorithmic}
{\small
\item
\texttt{INPUT:} $N$ -- the number of iterations \\
\hspace{11mm} $n$ -- the order of Picard iterations\\
\hspace{5mm} $x_0 \in X$, $\theta_0 \in \Theta$, ${{\mbf u}}_0 = x_0$, $\alpha  > 0, \sigma \in (0, 1]$, { $\varepsilon > 0$, $k = 0$}\\

 \texttt{REPEAT}\\
\qquad $\theta_{k+1} := \Pi_\Theta[ \theta_{k} - {\alpha \nabla_\theta L_{ls}(x_{k}, \theta_{k}, \mcl P^n_{x_{k}, \theta_{k}}{{\mbf u}}_{k})}]$ \\
\qquad $x_{k+1} :=  \Pi_X[x_{k} - {\alpha \nabla_{x} L_{ls}(x_{k}, \theta_{k}, \mcl P^n_{x_{k}, \theta_{k}}{{\mbf u}}_{k})}]$ \\
\qquad${{\mbf u}}_{k+1}(t)  := (1-\sigma) {{\mbf u}}_{k}(t) + \sigma(\mcl P_{x_{k+1}, \theta_{k+1}} {{\mbf u}}_{k})(t)$ \\
 \qquad $k = k + 1$,\\
{ \texttt{UNTIL} $\|\nabla_{x, \theta}L_{ls}(x_{k}, \theta_{k}, \mcl P^n_{x_{k}, \theta_{k}}{{\mbf u}}_{k})\| \leq \varepsilon$ or $k > N$.} \\
\texttt{OUTPUT:} $x_{N}, \theta_{N}, {{\mbf u}}_{N}$
}
\end{algorithmic}
\caption{Gradient-Contract-Multiplier Algorithm}  \label{alg3:modified_picard}
\end{algorithm}
\end{minipage}\vspace{-1mm}

\begin{minipage}[t]{0.5\textwidth}
\begin{algorithm}[H]
\begin{algorithmic}
{\small
\item
 \texttt{INPUT:} $N$ -- the number of iterations \\
 \hspace{11mm} $n$ -- the order of Picard iterations \\
 \hspace{11mm} $J$ -- the number of time intervals \\
 \hspace{11mm} $T > 0$, $\theta_0 \in \Theta$, $x_0 \in X^{\otimes J}$, ${\mbf u}_{0} \in C[0, T]^{\otimes J}$,  \\
 \hspace{11mm} $\lambda \geq 0$, $\alpha  > 0, \sigma \in (0, 1]$, { $\varepsilon > 0$}\\ 
 \texttt{INIT:} $\mbf u_j(t) = x_{0,j}$ for all $j \in \overline{1, J}$, { $k = 0$}\\
 \texttt{DEFINE:} $L_{\lambda, n}(x, \theta,\mbf u)$ as in Eqn.~\eqref{eqn:extended_loss}\\
 
\texttt{REPEAT}\\
\quad $\theta_{k+1} := \Pi_\Theta[ \theta_{k} - \alpha \nabla_\theta L_{\lambda, n} (x_{k}, \theta_{k},{\mbf u}_{k})]$ \\
\quad \texttt{For} $j$ from $1$ to $J$\\
\quad\quad $ {x}_{k+1, j} := \Pi_X[x_{k, j} - {\alpha \nabla_{x_j} L_{\lambda, n}(x_k, \theta_k, {\mbf u}_k )}]$ \\ 
\quad\quad $\mbf u_{k+1, j}  := (1-\sigma) \mbf u_{k, j} + \sigma(\mcl P_{(j-1)T, x_{k+1, j}, \theta_{k+1}} \mbf u_{k, j}) $ \\
\quad \texttt{EndFor}\\ 

 \quad $k = k + 1$,\\
{ \texttt{UNTIL} $\|\nabla_{x, \theta}L_{\lambda, n} (x_{k}, \theta_{k},{\mbf u}_{k})\| \leq \varepsilon$ or $k > N$.} \\
\texttt{OUTPUT:} $x_{N}, \theta_{N}, {{\mbf u}}_{N}$
}
\end{algorithmic}
\caption{Extended Grad-Contract-Multiplier Alg.}  \label{alg4:ext_picard}
\end{algorithm}
\end{minipage}\vspace{-3mm}

\end{figure} 
In this section we consider algorithms to solve the optimization problems in Sec.~\ref{sec:problem_statement}. For the unconstrained problem (Problems~\eqref{eqn:optimization_problem_unconstrained} and~\eqref{eqn:optimization_problem_unconstrained_approximated}), we focus on calculating the gradient of high-order Picard iterations. For the constrained optimization problem (Problems~\eqref{eqn:optimization_problem_constrained} and~\eqref{eqn:optimization_problem_constrained_approximated}), we focus on an iteration which eliminates the infinite-dimensional nature of the constraint and variables.

\subsection{A Gradient Descent Algorithm for the Unconstrained Formulation}\label{subsec:Alg1_grad_descent} As indicated in Subsec.~\ref{subsec:reformulation}, Optimization Problem~\eqref{eqn:optimization_problem_unconstrained} can be approximated using the Picard iteration as in~\eqref{eqn:optimization_problem_unconstrained_approximated}, where the Picard operator in this case is understood to act on the set of functions $\mbf u(t,x,\theta)$ and convergence of Picard iterates is then to the solution map. Unconstrained optimization problems of this form can be solved using gradient descent (combined with projection onto the convex set of allowable parameters, $\theta\in \Theta$ and initial conditions, $x \in X$). Specifically, let $\alpha > 0$ be a step size and recall that the projection operator ${\Pi}_X : \R^{n_x} \rightarrow \R^{n_x}$ is defined as\vspace{-3.5mm}
\begin{equation}\label{eqn:projection_operator}
    {\Pi}_X(z) := \underset{x \in X}{\texttt{argmin}} \|z - x\|^2_2\vspace{-3.5mm}
\end{equation}

We may now propose Alg.~\ref{alg1:gradient_algorithm}, where the gradient is computed using the n-th order Picard iterate approximation of the solution map. As stated in the following theorem, this gradient also converges to the gradient of the solution map -- implying that Alg.~\ref{alg1:gradient_algorithm} is actually an approximation of the gradient descent algorithm as applied the original formulation of the problem in Eqn.~\eqref{eqn:optimization_problem_unconstrained}.

\begin{thm}\label{thm:picard_sensitivity}
Let $\Gamma, X, \Theta$ be compact and $T$ be sufficiently small. Let $a = 2\sup_{x \in X} \norm{x}_2$ and $f \in C^1(\Gamma \times B_a \times\Theta)$ and $\nabla_x f, \nabla_\theta f$ are Lipschitz continuous functions. Then for all $t\in[0,T]$, $x \in X$, $\theta\in \Theta$ and for any  $\mbf u \in C([0, T])$ such that $\|\mbf u\|_\infty \leq a$ we have $\lim_{n\rightarrow \infty} \nabla_{x,\theta} ({\mcl P}_{ x, \theta}^n \mbf u)(t) = \nabla_{x, \theta} \phi(t , x, \theta)$,
where $\phi(t, x, \theta)$ is the solution map (Eqn.~\eqref{eqn:solutionmap}).
\end{thm}\vspace{-3.5mm}
\begin{pf}
See Prop.~\ref{prop:picard_sensitivity} in Sec.~\ref{sec:picard}. $\blacksquare$
\end{pf}\vspace{-3.5mm}
Convergence of the gradient descent in Alg.~\ref{alg1:gradient_algorithm} to the solution of Optimization Problem~\eqref{eqn:optimization_problem_unconstrained} as $n \rightarrow \infty$ is then guaranteed if sets $X, \Theta$ are convex and the loss function is strongly convex -- See, e.g., Thm. 4.32 in~\cite{bonnans2013perturbation}. Of course, in reality, for most parameter estimation problems the dependence of the loss function in Eqn.~\eqref{eqn:optimization_problem_unconstrained_approximated} on the system parameters and initial state is not convex, and hence Alg.~\ref{alg1:gradient_algorithm} is only guaranteed to converge to a set of local minima.
	
More significantly, convergence is only guaranteed in the limit $n \rightarrow \infty$. However, when $n$ becomes large, computation of the gradient $\nabla_{x,\theta} {\mcl P}_{x, \theta}^n$ becomes difficult. For this reason, we now consider an algorithm for solving the constrained form of Optimization Problem~\eqref{eqn:optimization_problem_constrained_approximated}.


\subsection{Gradient Contraction for Constrained Optimization}\label{subsec:Alg2_grad_contraction}

Now let us consider the constrained form of the optimization problem given in~\eqref{eqn:optimization_problem_constrained_approximated}. To account for the constraint, we use a two-step gradient contractive algorithm which includes both a gradient step on variables $x,\theta $ and a contraction step on variable $\mbf u$. In a standard gradient \textit{projection} algorithm, we would evaluate the gradient with respect to $x,\theta ,\mbf u$ and then project $\mbf u$ onto the feasible set. However, because $\mbf u$ is infinite-dimensional and the feasible set is a manifold on this infinite-dimensional space, such an approach is difficult. As an alternative, however, one might propose Alg.~\ref{alg2:picard}, which does not explicitly parameterize the variable $\mbf u$, but instead takes the gradient with respect to $x$ and $\theta$ for some fixed $\mbf u$ as \vspace{-3.5mm}
\begin{align*} 
\theta_{k+1} &= \Pi_\Theta[\theta_k - {\alpha \nabla_\theta L_{ls}(x_k, \theta_k, {{\mbf u}}_k)}] \\  x_{k+1} &= \Pi_X[x_k - {\alpha \nabla_{x} L_{ls}(x_k, \theta_k, {{\mbf u}}_k)}], \\[-8mm]
\end{align*} 
where note that, unlike in Alg.~\ref{alg1:gradient_algorithm}, there is no need to compute the gradient of a Picard iteration. Now, instead of parameterizing $\mbf u$ explicitly, for given $x,\theta $, we use the Picard iteration to update our variable $\mbf u$ as \vspace{-3.5mm}
\[
{{{\mbf u}}_{k+1}} = (1- \sigma) {{\mbf u}}_k + \sigma \mcl P_{ x_{k+1}, \theta_{k+1}} \mbf u_k\vspace{-3.5mm}
\]
 for some step size, $\sigma \in (0, 1]$. Then, because the Picard iteration is a contraction on sufficiently short time intervals, this step not only updates $\mbf u$, but moves it closer to the feasible set $\{\mbf u\;:\; \mbf u = {\mcl P}_{x_{k+1}, \theta_{k+1}} \mbf u\}$. 
 
Now, Lipschitz continuity and strong convexity of the loss function with respect to $x,\theta $ implies convergence of Alg.~\ref{alg2:picard} to a  fixed point $x^*_0$, $\theta^*_0$ and $\mbf u^*_0$. However, our failure to explicitly parameterize $\mbf u$ does have a cost. Specifically, while the fixed point will be feasible, this fixed point will not necessarily be optimal. To see this, consider the KKT conditions necessary for optimality of Problem~\eqref{eqn:optimization_problem_constrained_approximated}. 

\begin{prop}\label{prop:optimality_conditions}
Suppose $f, g \in C^1(\Gamma\times\R^{n_x}\times\Theta)$ and $\{x, \theta, \mbf u\} \in int(X) \times int(\Theta) \times C[0, T]$ is a stationary point of Optimization Problem~\eqref{eqn:optimization_problem_constrained}. Then \vspace{-5mm}
\begin{align} 
\nabla_{x, \theta} {L_{ls}}(x, \theta, \mbf u) + \sum_{i = 1}^{N_s}\nabla_{u_i} L_{ls}(x, \theta, \mbf u) \nabla_{x, \theta} \phi(t_i, x, \theta) = 0,\notag\\[-4mm]
\label{eqn:optimality_conditions}\\[-9mm]\notag
\end{align}
where $\phi(t, x, \theta)$ is the solution map of Differential Eqn.~\eqref{eqn:ODE} and \vspace{-5mm}
\begin{equation*} 
\nabla_{u_i} L_{ls}(x, \theta, \mbf u)\hspace{-1mm}:=\hspace{-1mm}\frac{1}{N_s}\hspace{-1mm} (y_i - g(t_i, \mbf u(t_i), \theta))^T \nabla_u g(t_i, \mbf u(t_i), \theta).\vspace{-3.5mm}
\end{equation*}
\end{prop}\vspace{-5mm}
\begin{pf}
The proof is based on equivalence of Problems~\eqref{eqn:optimization_problem_unconstrained_approximated} and~\eqref{eqn:optimization_problem_constrained} and application of the KKT conditions as formulated in~\cite{luenberger1997optimization}.$\blacksquare$
\end{pf}\vspace{-3.5mm}
Prop.\hspace{-0.5mm}~\ref{prop:optimality_conditions} implies the solution, $x^*\hspace{-0.5mm}, \theta^*\hspace{-0.5mm}, \mbf u^*$ of Problem\hspace{-0.5mm}~\eqref{eqn:optimization_problem_constrained} need satisfy Eqn.~\eqref{eqn:optimality_conditions}. However, (interior) fixed points of Alg.~\ref{alg2:picard} only satisfy $\nabla_{x, \theta} {L_{ls}}(x, \theta, \mbf u)=0$. To resolve this issue, in the following subsection, we approximate the multiplier $\nabla_{x, \theta} \phi(t_i, x, \theta)$ in Eqn.~\eqref{eqn:optimality_conditions} using Picard iteration and include this term in the updates of  $x$ and $\theta $.

\subsection{A Gradient Contraction Multiplier Algorithm for the Constrained Formulation} \label{subsec:Alg3_KKT_multiplier}
We now present Alg.~\ref{alg3:modified_picard}, which does not require the use of infinite-dimensional variables and constraints, but for which the fixed point satisfies an approximation of the KKT conditions for optimality defined in Eqn.~\eqref{eqn:optimality_conditions} of Prop.~\ref{prop:optimality_conditions}. Specifically, we require a fixed point $\{x,\theta ,\mbf u\}$ of Alg.~\ref{alg3:modified_picard} to satisfy\vspace{-3.5mm}
\begin{align}
&\nabla_{x,\theta} L_{ls}(x, \theta, \mcl P^n_{ x, \theta}\mbf u) = \nabla_{x,\theta} L_{ls}(x, \theta, \mbf u) \label{eqn:KKT_approximation}\\[-1mm]
&\quad + \sum\nolimits_{i = 1}^{N_s}\nabla_{u(t_i)} L_{ls}(x, \theta, \mbf u) \nabla_{x,\theta} (\mcl P^n_{ x, \theta}\mbf u)(t_i, x, \theta)=0\notag \\[-8mm]\notag
\end{align}
where recall $\lim_{n\rightarrow \infty} \nabla_{x,\theta} ({\mcl P}_{ x, \theta}^n \mbf u)(t) = \nabla_{x, \theta} \phi(t, x, \theta)$ as in Thm.~\ref{thm:picard_sensitivity}. Thus for any $n$, the optimality conditions of Prop.~\ref{prop:optimality_conditions} are satisfied in some approximate sense, where the accuracy of this approximation increases with $n$. As will be shown in Sec.~\ref{sec:main_algorithm} this allows us to prove convergence of the algorithm to optimality under suitable regularity and convexity conditions. As a practical matter, as will be shown in the numerical examples (Sec.~\ref{sec:examples}), Alg.~\ref{alg3:modified_picard} requires only $n=1$ or $n=2$ for highly accurate solutions.

\subsection{An Extended Gradient Contraction Multiplier Algorithm for the Constrained Formulation} \label{subsec:Alg4_KKT_multiplier_extended}
Alg.~\ref{alg2:picard} and Alg.~\ref{alg3:modified_picard} eliminated the need for explicit parametrization of the variable $\mbf u$ by initializing $\mbf u$ with some arbitrarily chosen initial guess and then using the Picard operator in lieu of a projection to enforce the equality constraint $\mcl P_{x, \theta} \mbf u - \mbf u=0$. However, this approach is premised on the assumption that the Picard operator is sufficiently contractive so that, after a sufficient number of iterations, this equality constraint will be satisfied. However, the time interval for which the Picard iteration is contractive is substantially smaller than the interval on which the data is typically sampled or for which the solution map is defined. As a final step, therefore, we partition the time-domain into disjoint intervals and apply the Picard iteration to each interval, using a regularization term to ensure that discontinuities between intervals are minimized. 

Specifically, for a given fixed interval of convergence, $t \in [0,T]$, we lift the space of initial conditions as $x\hspace{-0.5mm} \in \hspace{-0.5mm}X^{\otimes J}\hspace{-0.5mm}:=\hspace{-0.5mm}X \hspace{-0.5mm}\times \hspace{-0.5mm}\cdots\hspace{-0.5mm} \times\hspace{-0.5mm} X$ and the space of solutions as $\mbf u \in C[0,T]^{\otimes J}$ so that $\mbf u_{j}(0) = x_{j}$ and $\mbf u_{j}(T) = x_{j+1}$ for all $j \in \overline{1, J-1}$ and $\mbf u_j(t)$ represents the estimated solution on interval $t\in[(j-1)T, jT)$. Then the constrained optimization problem (Eqn.~\eqref{eqn:optimization_problem_constrained}) may now be formulated as \vspace{-3.5mm}
\begin{align} 
\hspace{-1mm}\min_{\substack{x \in X^{\otimes J},\, \theta \in \Theta\\ \mbf u \in C[0,T]^{\otimes J}}} &\;   \frac{1}{2N_s} \hspace{-1mm} \sum\nolimits_{i=1}^{N_s}\|y_i - g(t_i, \mbf u_{\lfloor t_i/T \rfloor + 1}(t_i\hspace{-0mm}\bmod{T}), \theta)\|_2^2\notag\\[-2mm] 
     \text{s.t.} \; & \mcl P_{(j-1)T, x_j, \theta} \mbf u_j = \mbf u_j \qquad \forall j ,\,  t\in [0,T] \label{eqn:optimization_problem_constrained_extended1} \\[-1mm]
      & \mbf u_{j}(T)  = x_{j+1} \qquad \forall j \in \overline{1, J-1}.\notag \\[-8mm] \notag
\end{align}
%
To see that Problems~\eqref{eqn:optimization_problem_constrained_approximated} and~\eqref{eqn:optimization_problem_constrained_extended1} are equivalent, we have the following.
\begin{lem} Let $\Gamma = [0, JT]$ and $\mbf u(jT) \in X$ for all $j \in \overline{1, J}$, then Optimization Problems~\eqref{eqn:optimization_problem_constrained_approximated} and~\eqref{eqn:optimization_problem_constrained_extended1} are equivalent.
\end{lem} \vspace{-3.5mm}
\begin{pf} Given a solution $\mbf u \in C[0, JT], x \in X, \theta \in \Theta$ to Problem~\eqref{eqn:optimization_problem_constrained_approximated} with objective value $\gamma$, let $x_j = \mbf u((j-1)T)$ and $\mbf u_j(t) = \mbf u(t + (j-1)T)$ for all $j\in\overline{1, J}$ and $t \in [0, T]$. Since $\mbf u$ is a solution of Problem~\eqref{eqn:optimization_problem_constrained_approximated}, we have $x_j = \mbf u((j-1)T) = (\mcl P_{0, x, \theta} \mbf u)((j-1)T)$ -- implying $\mbf u_{j}(T)  = x_{j+1}$. Next, since $\mcl P_{0, x, \theta} \mbf u = \mbf u$,  we have \vspace{-3.5mm}
\begin{align*}
& (\mcl P_{(j-1)T, x_j, \theta} \mbf u_j)(t) = (\mcl P_{(j-1)T, 0, \theta} \mbf u_j)(t) + x_j \\
&\quad= (\mcl P_{(j-1)T, 0, \theta} \mbf u_j)(t) + (\mcl P_{0, x, \theta} \mbf u)((j-1)T) \\
&\qquad =(\mcl P_{0, x, \theta} \mbf u)(t+ (j-1)T) = \mbf u_j(t)\\[-8mm]
\end{align*}
 for all $j \in \overline{1, J}$ and $t \in [0, T]$. Thus, $\mbf u_j, x_j, \theta$ is feasible for Problem~\eqref{eqn:optimization_problem_constrained_extended1} with objective value $\gamma$.

Conversely, given a solution $\mbf u_j,x$ to Problem~\eqref{eqn:optimization_problem_constrained_extended1} with objective value $\hat \gamma$, let $\hat{\mbf u}(t)=\mbf u_{\lfloor t/T \rfloor+1}(t \bmod{T})$ and $\hat x = x_1$. Then, since $\mbf u_j(T) = x_{j+1} = \mbf u_{j+1}(0)$, $\hat{\mbf u}$ is continuous.
 Now, suppose for some $j_0 \in \overline{1, J}$ and any $t \leq j_0 T$  we have $(\mcl P_{0, x_1, \theta} \hat{\mbf u})(t) = \hat{\mbf u}(t)$, which clearly holds for $j_0 = 1$.  Then, for all $t \in [j_0 T, (j_0+1)T]$ we have\vspace{-3.5mm}
\begin{align*}
&(\mcl P_{0, x_1, \theta} \hat{\mbf u})(t) = (\mcl P_{\lfloor t/T \rfloor \cdot T, x_{1}, \theta}  \mbf u_{\lfloor t/T \rfloor+1})(t \bmod T) \\ 
&\qquad\qquad\qquad\qquad\qquad- x_1 +(\mcl P_{0, x_{1}, \theta}  \hat{\mbf u})(j_0 T) \\ 
&\qquad=  (\mcl P_{\lfloor t/T \rfloor \cdot T, x_{{\lfloor t/T \rfloor}+1}, \theta}  \mbf u_{\lfloor t/T \rfloor+1})(t \bmod T) = \hat{\mbf u}(t). \\[-8mm]
\end{align*}

Hence, by induction we have $(\mcl P_{0, x_1, \theta} \hat{\mbf u})(t) = \hat{\mbf u}(t)$ for all $t \in [0, JT]$. We conclude that, $\hat{\mbf u}, \hat x, \theta$ is feasible for Problem~\eqref{eqn:optimization_problem_constrained_approximated} with objective value $\hat \gamma$. $\blacksquare$
 \end{pf}\vspace{-3.5mm}
 
To solve Optimization Problem~\eqref{eqn:optimization_problem_constrained_extended1}, we again use a Picard contractive step to avoid explicit parametrization of the variables $\mbf u_i$. However, this approach makes it difficult to directly enforce the linking constraints $\mbf u_{j}(T) = x_{j+1}$. Our approach, then is to move this constraint into the objective by using a penalty function as \vspace{-3.5mm}
\begin{align} \label{eqn:optimization_problem_constrained_extended}
 \hspace{-1mm}  \min_{\substack{x \in X^{\otimes J}, \theta \in \Theta\\ \mbf u \in C[0,T]^{\otimes J}}}\hspace{-0.5mm} &\; \frac{1}{2N_s} \hspace{-1mm} \sum\nolimits_{i=1}^{N_s} \hspace{-0.5mm} \|y_i\hspace{-0.5mm} -\hspace{-0.5mm} g(t_i, \mbf u_{\lfloor t_i/T \rfloor + 1}(t_i\hspace{-0mm}\bmod{T}), \theta)\|_2^2 \notag \\[-5mm]
    &\qquad\qquad + \lambda \sum\nolimits_{j=1}^{J-1}\|\mbf u_{j}(T) - x_{j+1}\|_2^2\\
     \text{s.t.} \;&\;  \mcl P_{(j-1)T, x_j, \theta} \mbf u_j  = \mbf u_j  \qquad \forall j \in \overline{1, J}, \notag\\[-8mm]\notag 
\end{align}
where $\lambda \geq 0$ is a regularization parameter. 


Having formulated the constrained optimization problem, we now propose an extended version of the algorithm discussed in Subsec.~\ref{subsec:Alg3_KKT_multiplier} where Picard iterations are introduced in the objective function so that the gradient of this objective approximates the KKT conditions in Eqn.~\eqref{eqn:optimality_conditions} -- similar to the approach described in Eqn.~\eqref{eqn:KKT_approximation}. Specifically, we now have objective function \vspace{-4mm}
\begin{align}\label{eqn:extended_loss}
& L_{\lambda, n}(x, \theta, \mbf u) \hspace{-1mm}:= \hspace{-1mm} \sum_{j=1}^{J-1} \lambda \|\mcl P^n_{(j-1)T, x_j, \theta} \mbf u_{j}(T) - x_{j+1}\|_2^2\\[-3mm]
&\hspace{-3mm} +\hspace{-1mm}\frac{1}{2N_s}\hspace{-1mm}  \sum_{i=1}^{N_s}\hspace{-1mm}\|y_i \hspace{-0.5mm}-\hspace{-0.5mm} g(t_i, (\mcl P^n_{(j-1)T, x_j, \theta} \mbf u_{\lfloor t_i/T \rfloor+ 1})(t_i\hspace{-0mm}\bmod{T}), \theta)\|_2^2 \notag  \\[-8mm] \notag
\end{align} 
for which we define Alg.~\ref{alg4:ext_picard}.

Although this approach allows us to approximate the solution of a parameter estimation problem on an arbitrary time interval, there are some technical challenges. First, the approach requires additional variables $x_{j}$ (initial conditions at time points $(j-1)T$) -- increasing the dimension of the optimization problem. Second, the use of a penalty function in lieu of the constraint $\mbf u_j(T) = x_{j+1}$ implies discontinuities unless $\lambda$ is chosen to be very large. However, increasing the value of $\lambda$ increases the gradient of the objective -- resulting in the need for smaller step sizes and hence slowing the convergence rate of the algorithm.  These numerical issues are analyzed in detail in Subsec.~\ref{subsec:VDP_reg}.
 
%
%

\section{Gradient-Contractive Algorithms}\label{sec:main_algorithm}
In this section, we propose sufficient conditions for convergence of the gradient-contractive algorithms presented in Sec.~\ref{sec:algorithms}. Specifically, we propose a general class of optimization problems which includes those posed in~\eqref{eqn:optimization_problem_constrained_approximated} and~\eqref{eqn:optimization_problem_constrained_extended}. Then, based on this generalized form of optimization problem, we define a generalized class of gradient-contractive algorithms. These two-step algorithms repeatedly update the optimization variables using both a gradient descent step and a contraction map (Algs.~\ref{alg2:picard},~\ref{alg3:modified_picard} and~\ref{alg4:ext_picard} are special cases of this approach). Finally, in Thm.~\ref{thm:gradient-contractive-map-is-contraction}, we provide sufficient conditions such that each step of the gradient-contractive algorithm is a contraction -- implying convergence.

To begin, suppose we are given a general optimization problem of the form\vspace{-3.5mm}
 %
\begin{equation} \label{eqn:general_optimization_problem}
   \min_{\substack{x \in X,  \mbf u \in U}} \; L(x, \mbf u)\quad 
     \text{s.t.} \;\;  \mcl P(x, \mbf u)  = \mbf u, \vspace{-3mm}
\end{equation}
where $ U \subset \mcl F(\R^{n_i}, \R^{n_u})$, $X \subset \R^{n_x}$, $\mcl P:X\times U \rightarrow U$  and  $L:X\times U \rightarrow \R$.
Now the gradient-contractive algorithm for step sizes $\alpha > 0$, $\sigma \in (0, 1]$ and order $n\in \N$ is defined by the sequence $(x_{k},\mbf u_{k})=\mcl T^k (x_{0},\mbf u_{0})$ where we say $(x_{k+1},\mbf u_{k+1})=\mcl T (x_{k},\mbf u_{k})$ if \vspace{-3.5mm}
\begin{align}\label{eqn:contractive_map}
x_{k+1} &=  \Pi_X\big[x_k - \alpha \nabla_x L(x_k, \mcl P^n(x_k, \mbf u_k))\big] \notag\\ 
\mbf u_{k+1} &=   (1 - \sigma) \mbf u_{k} + \sigma \mcl P(x_{k+1}, \mbf u_k), \\[-8mm] \notag
\end{align}
%
where recall $\Pi_X$ is the projection operator.

 Clearly, Algs.~\ref{alg3:modified_picard} and~\ref{alg4:ext_picard} are special cases of this generalized gradient contractive algorithm, where $n=0$ for Alg.~\ref{alg2:picard}, $\mcl P$ is the Picard operator, $L_i$ are as defined in Eqn.~\eqref{eqn:optimization_problem_unconstrained} and parameter variables $\theta$ are combined into $x$.  

Now, we provide conditions under which $\mcl T$ is a contraction, defined as follows.

 \begin{defn}[Contraction]\label{def:contraction}
     Given a metric space $H$, we say $\mcl T: H \rightarrow H$ is {a contraction} if there exists $q \in [0, 1)$ such that $\|\mcl T u_1 - \mcl T u_2\|_{H} \leq q \|u_1 - u_2\|_{H}$ for all $u_1, u_2 \in H$.
 \end{defn}
When proving contractivity, we use the norm $\norm{x,\mbf u}_H=\norm{x}_X+\norm{\mbf u}_U$.

The key property of contractive mappings is that, when iteratively applied to any $u\in H$, they converge to a \textit{fixed point}, $u^*=\lim_{k \rightarrow \infty} \mcl T^k u_0$, for which $\mcl T u^*=u^*$.

 \begin{thm}[Banach fixed-point theorem]\label{thm:banach-fixedpoint}
     Let $H$ be a complete metric space and $T:H \rightarrow H$ be a contraction. Then there exists unique fixed-point $u^* \in H$ such that $\mcl T u^* = u^*$. Furthermore,  $\lim_{k \rightarrow \infty} \mcl  T^k u_0 = u^*$ for any $u_0\in H$.
 \end{thm}

Thm.~\ref{thm:gradient-contractive-map-is-contraction} gives sufficient conditions, expressed in terms of the properties of $\mcl P$ and $L$, under which $\mcl T$ is a contraction and, which, therefore, implies that $\mcl T^k$ converges to a solution of Problem~\eqref{eqn:general_optimization_problem}.


\begin{thm}\label{thm:gradient-contractive-map-is-contraction}
Let  $U \subset \mcl U:=\mcl C_0(\R,\R^{n_u})$ and $X \subset \R^{n_x}$ be convex, closed and bounded sets with $r=\sup_{\mbf u\in U}\norm{u}_{\mcl U}$. Suppose that
\begin{enumerate}
 \item $\mcl P:X\times U \rightarrow U$ is differentiable in $x \in X$ and $\mcl P(x, \mbf  u)$ is a contraction in $\mbf u$, uniformly in $x \in X$. 

\item  $\nabla_x \mcl P(x, \mbf u)$ is  Lipschitz continuous on $(x,\mbf u) \in X \times U$ and there exist  $q < 1$, $K > 0, N\in \N$ such that for all $n \ge N$, $x_1,x_2\in X$ and $\mbf u_1, \mbf u_2 \in U$ \vspace{-3.5mm}
\begin{align}
\hspace{-10mm}\|\nabla_x \mcl P^n(x_1, \mbf u_1) - \nabla_x \mcl P^n(x_2, \mbf u_2)\| & \leq q^n \|\mbf u_1 - \mbf u_2\|_{\mcl U} \label{eqn:general_picard_is_lipschitz} \\ & \qquad + K \| x_1 -  x_2\|_{X}. \notag\\[-9mm]\notag 
\end{align} 
\item {  $L$ has form $L(x,\mbf u)=\sum_{i = 1}^{N_s} L_i(x, \mbf u(t_i))$ with $L_i(x,u) \in C^1(X \times B_r)$ and $\nabla L_i(x,u)$ Lipschitz on $(x,u) \in X \times B_r$ where $B_r:=\{u \in \R^{n_u}\;:\; \norm{u}\le r\}$.  
    \item $L(x,\mcl P^n(x,\mbf u))$ is strongly convex on $x\in X$, uniformly in $\mbf u \in U$ and $n \in \N$.}
\end{enumerate}
Then if $\mcl T$ is defined as in~\eqref{eqn:contractive_map}, there exist  $\alpha >0$, $\sigma \in (0, 1]$ and $n \in \N$ such that 
\begin{enumerate}
\item  [a)] $\mcl T: X \times  U\rightarrow X \times U$ and is contractive.
%
%
\item [b)]There exist $\nu < 1$ and $(x^*, \mbf u^*) \in X \times U$ such that for any $(x_0, \mbf u_0)\in X \times U$, $\lim_{k\rightarrow \infty} \mcl T^k (x_0, \mbf u_0) = (x^*, \mbf u^*)$ and \vspace{-3.5mm}
    \[
\hspace{-7mm}\|\mcl T^k \hspace{-0.5mm}(x_0, \mbf u_0)\hspace{-0.1mm} -\hspace{-0.1mm} (x^*\hspace{-0.5mm}, \mbf u^*\hspace{-0.5mm})\|_{X\times \mcl U} \hspace{-0.5mm} \leq \hspace{-0.5mm}\nu^k(\|x_0 - x^*\hspace{-0.5mm}\|_X + \|\mbf u_0 - \mbf u^*\hspace{-0.5mm}\|_{\mcl U}).\vspace{-3.5mm}
\]
\end{enumerate}
\end{thm}
\begin{pf}{ 
The proof consists of three parts. First, we show that by condition 2), $L(x,\mbf u)$  is Lipschitz on $(x,\mbf u) \in X \times U$ with Lipschitz constant decreasing as $n$ increases. Second, we give conditions on $\alpha,\sigma$ and $n$ for the iterate map, $\mcl T$, to be contractive. Next, we prove existence of $\alpha, \sigma, n \in \N$ which satisfy these conditions. The theorem statement then follows directly from contractivity of $\mcl T$ by Thm.~\ref{thm:banach-fixedpoint}. 
 
\textbf{Variation of $\, L(x,\mbf u)$ with $\mbf u$.}  By Condition 2), we have that $\mcl P(x, \mbf u)$ is Lipschitz continuous in $x$ and there exists $q_1 < 1$ and $K_1$ such that  \vspace{-3.5mm}
\[
\|\mcl P(x_1, \mbf u_1) - \mcl P(x_2, \mbf u_2)\|_{\mcl U} \leq q_{1} \|\mbf u_1 - \mbf u_2\|_{\mcl U}  + K_{1} \|x_1 - x_2\|_X, \vspace{-3.5mm} \]
for any $x_1,x_2 \in X$ and $\mbf u_1,\mbf u_2 \in U$.

Now, using $q, K$ as in Condition 2), define $q_{\mcl P}=\max\{q_1,q\}<1$ and $K_{\mcl P}=\max\big\{\frac{K_1}{1 - q_1},K\big\}$. Then from
 Eqn.~\eqref{eqn:general_picard_is_lipschitz},  for all  $n \ge N$, $x_1, x_2 \in X$ and $\mbf u_1, \mbf u_2 \in U$, \vspace{-3.5mm}
 \begin{align*}
\|\mcl P^n(x_1, \mbf u_1) - \mcl P^n(x_2, \mbf u_2)\|_{\mcl U} &\leq q^n_{\mcl P} \|\mbf u_1 - \mbf u_2\|_{\mcl U} \\
& \quad  + K_{\mcl P} \|x_1 - x_2\|_X \\
\|\nabla_x \mcl P^n(x_1, \mbf u_1) - \nabla_x \mcl P^n(x_2, \mbf u_2)\|_{\mcl U} &\leq q^n_{\mcl P} \|\mbf u_1 - \mbf u_2\|_{\mcl U} \\
& \quad  + K_{\mcl P} \|x_1 - x_2\|_X. \\[-8mm]
\end{align*}
The latter inequality also implies $\|\nabla_x \mcl P(x, \mbf u)\|\leq  K_{\mcl P}$.

Now, since $\nabla_x L_i$ are Lipschitz continuous, by using the triangle inequality and summing over all Lipschitz constants, we may find $K_L \geq 0$ and $K_{n, \mbf u}$  such that \vspace{-3.5mm} 
\begin{align*}
&\|\nabla_x L(x_1, \mcl P^n(x_1, \mbf u_1)) -\nabla_x  L(x_2, \mcl P^n(x_2, \mbf u_2))\|_2 \\
&\qquad\qquad \leq K_{n,\mbf u} \|\mbf u_1 - \mbf u_2\|_{\mcl U} + K_{L} \|x_1 - x_2\|_X \\[-8mm]\notag
\end{align*}
 for all $n \in \N$, $x_1, x_2 \in X$ and $\mbf u_1, \mbf u_2 \in U$ where $q_{\mcl P}<1$ implies $\lim_{n\rightarrow \infty} K_{n,\mbf u} = 0$ --- a full proof of this statement can be found in Lem.~16 in Appendix A.


\textbf{Sufficient Conditions for $\mcl T$ to be Contractive.} Now consider the iterate map, $\mcl T$.
Denote the first step in Eqn.~\eqref{eqn:contractive_map} as $x_{k+1}=\mcl D_n(x_k,\mbf u_k)$ where  ${\mcl D_n}(x, \mbf u) := \Pi_X[x - \alpha \nabla_x L(x, \mcl P^n(x, \mbf u))]$. 
%

Since $X$ is convex, the projection $\Pi_X$ is nonexpansive and hence Lipschitz continuous with factor 1. Furthermore, by Condition 4),  $ L(x, \mcl P^n(x, \mbf u))$ is strongly convex on $X$ uniformly in $\mbf u$ and $n$. Thus, there exists $\mu > 0$ such that $L(x, \mcl P^n(x, \mbf u))$ is strongly convex with modulus $\mu >0$ for any  $\mbf u \in U$ and $n \in \N$. Recalling, $\nabla_x L(x, \mcl P^n(x,\mbf u)$ is Lipschitz continuous with bound $K_L$ as established above. Thus, we may directly apply Thm.~1 in Sec.~7.2.1 in~\cite{polyak1987introduction} to show that for $q_D(\alpha) = \max\{|1 - \alpha \mu|, |1 - \alpha K_L|\}$,\vspace{-3.5mm}
\begin{align*}
\|{\mcl D_n}(x_1, \hspace{-0.5mm}\mbf u_1)\hspace{-0.5mm}-\hspace{-0.5mm} {\mcl D_n}(x_2,\hspace{-0.5mm} \mbf u_2)\| &\leq q_{D}(\alpha) \|x_1 \hspace{-0.5mm}-\hspace{-0.5mm} x_2\|\hspace{-0.5mm}\\
& + \hspace{-0.5mm}\alpha K_{n, \mbf u} \|\mbf u_1\hspace{-0.5mm} -\hspace{-0.5mm} \mbf u_2\|,\\[-8mm]
\end{align*}
 for all $x_1, x_2 \in X$ and $\mbf u_1, \mbf u_2 \in U$.
  
The map $\mcl T$ defined as in~\eqref{eqn:contractive_map} can be expressed as\vspace{-3.5mm}
\[
\mcl T(x, \mbf u) = \bmat{\mcl D_n(x, \mbf u) \\[-1mm]  (1-\sigma)\mbf u + \sigma\mcl P(\mcl D_n(x, \mbf u), \mbf u) }, \vspace{-3.5mm}
\]
where $\mcl T:X\times U\rightarrow X\times U$. 
%
Finally, to show that $\mcl T$ is Lipschitz continuous jointly in $x,\mbf u$, we apply the following sequence of relatively straightforward inequalities. First, \vspace{-3.5mm}
\begin{align*}
& \|{\mcl T}(x_1, \mbf u_1)\hspace{-1mm}-\hspace{-1mm}{\mcl T}(x_2, \mbf u_1)\|_{X\times \mcl U}  = \|{\mcl D_n}(x_1, \mbf u_1)\hspace{-1mm}-\hspace{-1mm}{\mcl D_n}(x_2, \mbf u_1)\|_X \\
&\quad+ \|\sigma{\mcl P}({\mcl D_n}(x_1, \mbf u_1), \mbf u_1) - \sigma{\mcl P}({\mcl D_n}(x_2, \mbf u_1), \mbf u_1)\|_{\mcl U} \\ 
&\quad \leq q_{\mcl D}(\alpha) \|x_1\hspace{-1mm} -\hspace{-1mm} x_2\|_X + \sigma K_{\mcl P} \|{\mcl D_n}(x_1, \mbf u_1)\hspace{-1mm} -\hspace{-1mm} {\mcl D_n}(x_2, \mbf u_1)\|_X \\
&\quad \leq q_{\mcl D}(\alpha) \|x_1 - x_2\|_X + \sigma q_{\mcl D}(\alpha) K_{\mcl P} \|x_1 - x_2\|_X.  \\[-8mm]\notag
\end{align*}
for all $x_1, x_2 \in X$ and $\mbf u_1 \in \mcl U$. Next, \vspace{-3.5mm}
\begin{align*}
& \|{\mcl T}(x_2, \mbf u_1) - {\mcl T}(x_2, \mbf u_2)\|_{X\times \mcl U }  \\ 
& \;\leq \|{\mcl D_n}(x_2, \mbf u_1) - {\mcl D_n}(x_2, \mbf u_2)\|_{\mcl U} + (1-\sigma)\|(\mbf u_1 - \mbf u_2)\|_{\mcl U} \\
&\qquad\qquad+ \sigma\|{\mcl P}({\mcl D_n}(x_2, \mbf u_1), \mbf u_1) - {\mcl P}({\mcl D_n}(x_2, \mbf u_2), \mbf u_2)\|_{\mcl U} \\
&\; \leq (\alpha K_{n, \mbf u} + (1-\sigma) + \sigma q_{\mcl P})  \|\mbf u_1 - \mbf u_2\|_{\mcl U} \\
&\qquad\qquad + \sigma K_{\mcl P} \|{\mcl D_n}(x_2, \mbf u_1) - {\mcl D_n}(x_2, \mbf u_2)\|_X  \\
&\;\leq (1-\sigma + \sigma q_{\mcl P} + \alpha K_{n, \mbf u} +  \alpha \sigma K_{\mcl P} K_{n, \mbf u}) \|\mbf u_1 - \mbf u_2\|_{\mcl U}  \\[-8mm]\notag
\end{align*}
for all $x_2 \in X$ and $\mbf u_1, \mbf u_2 \in \mcl U$. And, finally, \vspace{-3.5mm}
\begin{align*}
&\|{\mcl T}(x_1, \mbf u_1) \hspace{-1mm}-\hspace{-1mm} {\mcl T}(x_2, \mbf u_2)\|_{X\times \mcl U} \hspace{-1mm}\\
\hspace{-1mm}&\leq\hspace{-1mm} \|{\mcl T}(x_1,\hspace{-0.5mm} \mbf u_1)\hspace{-1mm} - \hspace{-1mm}{\mcl T}(x_2,\hspace{-0.5mm} \mbf u_1)\|_{X\hspace{-1mm}\times \mcl U}\hspace{-1mm} +\hspace{-1mm} \|{\mcl T}(x_2, \hspace{-0.5mm}\mbf u_1)\hspace{-1mm} - \hspace{-1mm}{\mcl T}(x_2,\hspace{-0.5mm} \mbf u_2)\|_{X\hspace{-1mm}\times \mcl U} \\
&\leq\hspace{-1mm} q_{\mcl D}(\alpha)(\sigma K_{\mcl P} + 1)\|x_1 - x_2\|_X \\ 
&\qquad\quad + (1- \sigma + \sigma q_{\mcl P} + \alpha K_{n, \mbf u}(\sigma K_{\mcl P} + 1))\|\mbf u_1 - \mbf u_2\|_{\mcl U}  \\[-8mm]
\end{align*}
for all $x_1, x_2 \in X$ and $\mbf u_1, \mbf u_2 \in \mcl U$.
%
%
Therefore, we find that $\mcl T$ is a contraction if \vspace{-3.5mm}
\begin{equation*} 
(1- \sigma) + \sigma q_{\mcl P} + \alpha K_{n, \mbf u}(\sigma K_{\mcl P} + 1)< 1 \quad q_{\mcl D}(\alpha)(\sigma K_{\mcl P} + 1)  < 1.\vspace{-3.5mm}
\end{equation*}

\textbf{Existence of $\alpha, \sigma$ for Contractivity of $\mcl T$.} 

As our final step, we establish the existence of $\alpha, \sigma, n$ which satisfy the inequalities defined above and obtain an associated expression for the contractivity of $\mcl T$. First, we select $\alpha = \frac{2}{\mu + K_L}$, so that $q_{\mcl D}(\alpha) = \frac{K_L - \mu}{K_L + \mu}\hspace{-1mm} <\hspace{-1mm} 1$.\\  Next, since $\lim_{n \rightarrow \infty}K_{n,\mbf u} = 0$, there exists $n \geq 0$ such that $ K_{n,\mbf u} < \frac{(K_L - \mu)(1 - q_{\mcl P})}{2} \min\big\{\frac{2 \mu}{ K_{\mcl P}(K_L - \mu)}, 1 \big\}$ which implies $\frac{2 K_{n,\mbf u}}{(K_L - \mu)(1 - q_{\mcl P})}  < \min\big\{\frac{2 \mu}{ K_{\mcl P}(K_L - \mu)}, 1 \big\}$.


 

Now, choosing any $\sigma \hspace{-1mm} \in \hspace{-1mm}\big(\hspace{-0.15mm}\frac{2 K_{n,\mbf u}}{(K_L\hspace{-0.15mm} -\hspace{-0.15mm} \mu)(1\hspace{-0.15mm} -\hspace{-0.15mm} q_{\mcl P})},\min\hspace{-0.5mm}\big\{\hspace{-0.15mm}\frac{2 \mu}{ K_{\mcl P}(K_L\hspace{-0.15mm} -\hspace{-0.15mm} \mu)}, \hspace{-0.5mm}1\hspace{-0.5mm} \big\}\hspace{-0.5mm}\big)$, we define\vspace{-4.5mm}
\[
\nu_1 := q_{\mcl D}(\alpha)(\sigma K_{\mcl P} + 1) <  \frac{K_L\hspace{-1mm}- \hspace{-1mm}\mu}{K_L\hspace{-1mm} +\hspace{-1mm} \mu}\bigg(\hspace{-1mm}\frac{2 \mu}{ K_{\mcl P}(K_L\hspace{-1mm}-\hspace{-1mm}\mu)} K_{\mcl P} + 1\bigg) = 1,\vspace{-5mm}
 \]
and\vspace{-3.5mm}
\begin{align*}
& \nu_2 := (1- \sigma) + \sigma q_{\mcl P} + \alpha K_{n, \mbf u}(\sigma K_{\mcl P} + 1)  \hspace{-1mm} \\[-1mm]
& \qquad< \hspace{-1mm}  1 - \sigma(1  -  q_{\mcl P}) + \frac{2K_{n, \mbf u}}{K_L - \mu} \\[-1.5mm]
& \qquad < 1 - \frac{2 K_{n,\mbf u}}{(K_L - \mu)(1 - q_{\mcl P})}(1 - q_{\mcl P}) + \frac{2K_{n, \mbf u}}{K_L - \mu} = 1,  \\[-8mm] 
\end{align*}
then for $\nu:=\max\{\nu_1,\nu_2\}<1$, we have from the final inequality in part 2 of the proof that\vspace{-3.5mm}
\begin{align*}
&\|{\mcl T}(x_1, \mbf u_1) \hspace{-1mm}-\hspace{-1mm} {\mcl T}(x_2, \mbf u_2)\| \leq\hspace{-1mm} \nu\left(\|x_1 - x_2\|  + \|\mbf u_1 - \mbf u_2\|\right) \\[-8mm]
\end{align*}
for all $x_1, x_2 \in X$ and $\mbf u_1, \mbf u_2 \in \mcl U$.

We conclude that there exist $\alpha, \sigma, n$ such that the map $\mcl T$ is contractive with factor $\nu<1$, verifying statements a) and b). $\blacksquare$\vspace{-2.0mm}
%
%
%
%
}
\end{pf}\vspace{-2mm} 

 
 \textbf{Remark 1.} The gradient-contractive approach of Thm.~\ref{thm:gradient-contractive-map-is-contraction} may be extended to other gradient-based algorithms (defined by alternative maps $\mcl D_n$) as long as $\mcl D_n$ is contractive in $x$ and Lipschitz in $\mbf u$.

%

To apply the results of Thm.~\ref{thm:gradient-contractive-map-is-contraction} to the parameter estimation problem, we define the mapping $\mcl P$ as a Picard iteration. In the following section, we will recall properties of the Picard iteration. Then, in Sec.~\ref{sec:conditions_for_Algs}, we will give conditions on the parameterized vector field, $f$, under which the proposed modified and extended gradient contractive algorithms converge.


\section{Contraction and Lipschitz continuity of the Picard Operator}\label{sec:picard}

In Sec.~\ref{sec:main_algorithm}, we proposed a generalized form of gradient-contractive algorithm (Eqn.~\eqref{eqn:contractive_map}) for solving the class of optimization problems defined in Eqn.~\eqref{eqn:general_optimization_problem}, which are defined in terms of a Lipschitz separable objective, $L(x, \mbf u)$ and an equality constraint of the form $\mcl P(x,\mbf u)=\mbf u$. Furthermore, in Thm.~\ref{thm:gradient-contractive-map-is-contraction}, we provided conditions on $\mcl P$ and $L$ under which the gradient-contractive algorithm converges. In this section, we return to the original parameter estimation problem, as formulated in Optimization Problems~\eqref{eqn:optimization_problem_constrained_approximated} and~\eqref{eqn:optimization_problem_constrained_extended}, where $\mcl P$ is now the Picard operator. In this case, we show that Algs.~\ref{alg3:modified_picard} and~\ref{alg4:ext_picard} are special cases of the generalized gradient-contractive algorithm and that the Picard operator, $\mcl P_{t_0, x, \theta}$, satisfies the conditions of Thm.~\ref{thm:gradient-contractive-map-is-contraction}. Specifically, in Thm.~\ref{thm:picard_contraction} and Lem.~\ref{lem:picard_is_lipschitz_in_parameters} (Subsec.~\ref{subsec:condition1}), we show that $\mcl P_{t_0, x, \theta}$ satisfies Condition 1) of Thm.~\ref{thm:gradient-contractive-map-is-contraction} and in Lem.~\ref{lem:picard_sensitivity_is_continuous} (Subsec.~\ref{subsec:condition2}), we show that $\mcl P_{t_0, x, \theta}$ satisfies Condition 2). Next, in Sec.~\ref{sec:conditions_for_Algs}, we give conditions under which $L_{\lambda,n}$ (as defined in Eqn.~\eqref{eqn:extended_loss}) likewise satisfies the conditions of Thm.~\ref{thm:gradient-contractive-map-is-contraction}.
%

\subsection{Picard Iterations as a Contraction Map}\label{subsec:condition1}
In this subsection, we examine the Picard operator and show it satisfies Condition 1) of Thm.~\ref{thm:gradient-contractive-map-is-contraction}. Specifically, recall from Definition~\ref{def:picard_operator} that for any $t_0,x,\theta$ and parameterized vector field, $f$, the Picard operator is defined as \vspace{-3.5mm}
 \begin{equation*}
 (\mcl P_{t_0, x, \theta} \mbf u)(t) := x+ \int_{0}^{t} f(s + t_0, \mbf u(s), \theta) ds.\vspace{-3.5mm}
 \end{equation*} 

The contractive properties of the Picard operator and its gradient are relatively well-established on a sufficiently short time interval. However, in our case, we require these contractive properties to be uniform on domains $x\in X$ and $\theta \in \Theta$ which results in a slight variation on the standard result.
%
%
%

\begin{thm}[Picard-Lindel\"of Theorem]\label{thm:picard_contraction}
Suppose $\Gamma \subset \R$, $X\subset\R^{n_x}$, $\Theta\subset\R^{n_{\theta}}$ are compact and $f \in C(\Gamma\times\R^{n_x}  \times\Theta)$ is locally Lipschitz. Let $t_0 \in \Gamma$, $a=2\sup_{x \in X}\norm{x}_X$ and $T < \min\{\frac{1}{K_x}, \frac{a}{2M} \}$ where \vspace{-3.5mm}
\begin{equation*}
K_x = \sup_{\substack{ x, y \in B_{a},\, t\in\Gamma, \theta\in\Theta}} \frac{\|f(t, x, \theta) - f(t, y, \theta)\|_2}{\|x - y\|_X} \vspace{-3.5mm}
\end{equation*}
\begin{equation*}
M = \sup_{\substack{x \in B_{a}, t\in\Gamma, \theta\in \Theta}} \|f(t, x, \theta)\|_2. \vspace{-2mm}
\end{equation*}
Define $C_a[0, T]  : = \{\mbf u \in C([0,T]) \; | \;  \|\mbf u\|_{\infty} \leq a\}.$
Then for all $x \in X$ and $\theta \in \Theta$, ${\mcl P}_{t_0, x, \theta}:C_a[0, T] \rightarrow C_a[0, T]$ and there exists $\mbf u^* \in C_a[0, T]$ such that for any $\mbf u_0 \in C_a[0, T]$,  $\lim\limits_{k \rightarrow \infty}{\mcl P}_{t_0, x, \theta}^k \mbf u_0=\mbf u^*$.
Furthermore, $\|{\mcl P}_{t_0, x, \theta}\mbf u - {\mcl P}_{t_0, x, \theta} \mbf v \|_{\infty} \leq TK_x \|\mbf u - \mbf v\|_{\infty}$
for any $\mbf u,\mbf v \in C_a[0, T]$ 
%
\end{thm}\vspace{-3.5mm}

 \begin{pf}
Let $x \in X$ and $\theta \in \Theta$. We first show that ${\mcl P}_{t_0, x, \theta}:C_a[0, T] \rightarrow C_a[0, T]$.
 Recall $({\mcl P}_{t_0, x, \theta}\mbf u)(t) := x + \int_{0}^t f(s + t_0, \mbf u(s), \theta) ds$. Suppose $\mbf u \in C_a[0, T]$. Since $f$ is continuous, by, e.g.~\cite{rudin1964principles}, ${\mcl P}_{t_0, x, \theta}\mbf u$ is continuous.
 
  Moreover,\vspace{-5.5mm}
 \begin{align*}
 \|{\mcl P}_{t_0, x, \theta}\mbf u\|_{\infty} &= \sup_{t\in [0, T]} \hspace{-1mm} \big\|x + \int_{0}^t f(s + t_0, \mbf u(s), \theta) ds \big\|_2 \\[-4mm]
 &\leq \|x\|_2 + \sup_{ t \in [0, T]} \big\| \int_{0}^t f(s + t_0, \mbf u(s), \theta) ds \big\|_2 \\[-2mm]
 &\quad \leq \frac{a}{2} + T \hspace{-5mm} \sup_{\substack{ x \in B_{a} , \theta\in \Theta , s\in \Gamma}} \hspace{-5mm}\|f(s, x, \theta)\|_2 \leq \frac{a}{2} + TM.  \\[-9mm]\notag
 \end{align*}
Since $T\hspace{-1mm}\leq \hspace{-1mm}\frac{a}{2M}$, then $\|{\mcl P}_{t_0, x, \theta}\mbf u\|_{\infty} \hspace{-0.5mm}\leq\hspace{-0.5mm} a$ hence ${\mcl P}_{t_0, x, \theta} \mbf u\hspace{-1mm} \in\hspace{-1mm} C_a[0, T]$. 

The second part of this proof is to show the contraction property. For all $\mbf u, \mbf v \in C_a[0, T]$ we have \vspace{-3.5mm}
\begin{align*}
&\|{\mcl P}_{t_0, x, \theta}\mbf u - {\mcl P}_{t_0, x, \theta}\mbf v\|_{\infty}   \\
& = \sup_{t\in [0, T]}  \bigg\|\int_{0}^t f(s + t_0, \mbf u(s), \theta) - f(s + t_0, \mbf v(s), \theta) ds\bigg\|_2. \\[-9mm]\notag 
\end{align*}
Next, for all $\mbf u, \mbf v \in C_a[0, T]$ we have\vspace{-3.5mm}
\begin{align*}   
& \|{\mcl P}_{t_0, x, \theta}\mbf u - {\mcl P}_{t_0, x, \theta}\mbf v\|_{\infty}  \\[-1.5mm]
&\quad\leq \sup_{t \in \Gamma}  \int_{0}^t \|f(s + t_0, \mbf u(s), \theta) - f(s + t_0, \mbf v(s), \theta)\|_2 ds \\[-2mm]
&\qquad\qquad\qquad \leq T K_{x} \sup_{s\in \Gamma}   \|\mbf u(s) - \mbf v(s)\|_2 . \\[-9mm]\notag
\end{align*} 
Since $TK_{x} < 1$, we have ${\mcl P}_{t_0, x, \theta}$ is a contraction map. Since $C_a[0, T]$ is a closed and bounded  subset of the Banach space $C([0, T])$  and hence 
by the fixed point theorem, for any $\mbf u_0\in C_a[0, T]$, we have $\lim_{k \rightarrow \infty}{\mcl P}_{t_0, x, \theta}^k \mbf u_0=\mbf u^*\in C_a[0, T]$.$\blacksquare$
 \end{pf}\vspace{-4mm} 
Thm.~\ref{thm:picard_contraction} shows that $\mcl P_{t_0, x, \theta}$ satisfies Condition 1) of Thm.~\ref{thm:gradient-contractive-map-is-contraction}. 
Before moving on to prove that Condition 2) is satisfied (in Subsec.~\ref{subsec:condition2}), we state two further results which will be useful in establishing this result.
 
First, we note that, trivially, the solution map, $\phi$, itself is a fixed point of the Picard iteration, which is formally stated  as follows.
\begin{cor}\label{cor:picard_projects_onto_set_of_feasible_solutions}
Suppose the conditions of Thm.~\ref{thm:picard_contraction} are satisfied. Then for any $x\in X$, $\theta \in \Theta$  and for all $\mbf u \in C_a[0, T]$ \vspace{-3.5mm}
\begin{align*}
& \sup_{t \in [0, T]}\|(\mcl P_{0,x, \theta} \mbf u)(t) - \phi(t, x, \theta)\|_2 \\[-3mm]
& \quad\qquad \leq TK_{x} \sup_{t \in [0, T]}\| \mbf u(t) - \phi( t, x, \theta)\|_2,  \\[-8.5mm]\notag
\end{align*} 
where $\phi$ is the solution map of the ODE defined by $f$.
\end{cor}

Second, we establish a uniform Lipschitz continuity bound on $\mcl P_{t_0, x, \theta}$.


\begin{lem}[Lipschitz Continuity] ~\label{lem:picard_is_lipschitz_in_parameters}
Suppose $\Gamma, X, \Theta$ satisfy the conditions of Thm.~\ref{thm:picard_contraction} with $t_0$, $a$, $T$, $K_x$, $M$, $C_a[0, T]$ as defined therein.
Let\vspace{-3.5mm}
 $$
K_{\theta} = \sup_{x\in B_{a}, \; t\in\Gamma, \; \theta_1, \theta_2\in\Theta} \frac{\|f(t, x, \theta_1) - f(t, x, \theta_2)\|_2}{\|\theta_1 - \theta_2\|_2}. \vspace{-3.5mm}
$$
then for any $\mbf u \in C_a[0, T]$ and $n\in \N$, $\mbf v(t, x, \theta) := ({\mcl P}_{t_0, x, \theta}^n \mbf u)(t)$ is Lipschitz continuous in $x\in X$ and $\theta\in \Theta$  with bounds on the Lipschitz constants given by $\frac{1}{1 - TK_{x}}$ and $\frac{TK_{\theta}}{1 - TK_{x}}$, respectively.
   
\end{lem}\vspace{-4mm}
\begin{pf}
Define the sequence $\mbf v_0=\mbf u$, $\mbf v_{i}={\mcl P}_{t_0, x, \theta}^i \mbf u$.  Clearly  $\mbf v_0=\mbf u\in C_a[0, T]$ satisfies the Lipschitz bound since it is not a function of $x,\theta$. Now suppose that $\mbf v_i$ satisfies the Lipschitz bound. 
For all $x_1, x_2 \in X$, $\theta_1, \theta_2 \in \Theta$ we have \vspace{-3.5mm}
\begin{align*}
	& \sup_{t\in[0, T]}\|\mbf v_{i+1}(t,   x_1, \theta_1) -  \mbf v_{i+1}(t,  x_2, \theta_2)\|_2   \\[-3mm]
	 & =  \sup_{t\in[0, T]} \big\|x_1 - x_2 +  \int_{0}^{t} f(t_0 + s, \mbf v_{i}( s,  x_1, \theta_1), \theta_1) \\[-3mm]
	 & \qquad\qquad\qquad\qquad\qquad - f(t_0 + s, \mbf v_{i}(s,    x_2, \theta_2), \theta_2)ds\big\|_2 \\ 
	& \leq  TK_{x}  \sup_{t\in[0, T]}  \|\mbf v_{i}(t,   x_1, \theta_1) -\mbf v_{i}(t,   x_2, \theta_2)\| \\[-3mm]
	 & \qquad\qquad\qquad\qquad + \|x_1 - x_2\|_2 + TK_{\theta}\|\theta_1 - \theta_2\|_2.  \\[-8mm]\notag
\end{align*}

Since  $TK_x<1$, by induction we have \vspace{-3.5mm}
\begin{align*}
  &   \sup_{t\in[0, T]}\|\mbf v_n(t,   x_1, \theta_1) -  \mbf v_n(t,  x_2, \theta_2)\|_2   \\[-6mm]
  & \qquad\qquad \leq (\|x_1 - x_2\|_2 + TK_{\theta}\|\theta_1 - \theta_2\|_2)\sum_{k = 0}^{n-1} (TK_{x})^k \\[-2mm]
  & \qquad\qquad \leq \frac{\|x_1 - x_2\|_2 + TK_{\theta}\|\theta_1 - \theta_2\|_2}{1 - TK_{x}}  \\[-8mm]\notag
\end{align*}
for all $x_1, x_2 \in X$, $\theta_1, \theta_2 \in \Theta$ and $\mbf u \in C_a[0, T]$. Thus, $\mcl P_{t_0, x, \theta}^n \mbf u$ is Lipschitz, uniformly in $t_0,x,\theta \hspace{-.5mm}\in\hspace{-.5mm}\Gamma \hspace{-.5mm}\times\hspace{-0.5mm} X \hspace{-.5mm}\times\hspace{-.5mm} \Theta $, and $n\hspace{-0.5mm}\in \hspace{-0.5mm}\N$. $\blacksquare$  
\end{pf}\vspace{-3.5mm} 
Using these results, we now show that Condition 2) of Thm.~\ref{thm:gradient-contractive-map-is-contraction} is satisfied.
 


\subsection{Lipschitz continuity of Gradients of Picard Iterations}\label{subsec:condition2}

In this subsection, we show that $\mcl P$ satisfies Condition 2) of Thm.~\ref{thm:gradient-contractive-map-is-contraction}. First, 
we use an extended version of the result in~\cite{rudin1964principles} (Thm. 6.20) to show that if the parameterized vector field is sufficiently smooth, the Picard operator, $\mcl P$ is differentiable in $x,\theta$. 



 
 \begin{lem}[Rudin~\cite{rudin1964principles}]\label{lem:picard_is_continuous}
   Let $\Gamma, X$ and $\Theta$ be compact, and $f\in C^m(\Gamma\times \R^{n_x} \times \Theta)$ for $m \in \N$. If $\mbf u\in C^m(\Gamma\times X\times \Theta)$ and $\mbf v(t, x, \theta) := (\mcl P_{t_0, x, \theta}\mbf u)(t)$, then  $\mbf v \in C^m(\Gamma\times X\times \Theta)$.
\end{lem}  
Iteratively applying Lem.~\ref{lem:picard_is_continuous}, we have that $\mcl P^n$ is continuously differentiable with respect to $t_0,x,\theta \in \Gamma \times X\times \Theta$. If $f \in C^1$, we may thus define the gradients of an $n$-th order Picard operator as follows. \vspace{-3.5mm}
\begin{equation}
S_n(t, x, \theta) = \bmat{\nabla_x ({\mcl P}_{t_0, x, \theta}^n \mbf u)(t) & \nabla_\theta ({\mcl P}_{t_0, x, \theta}^n \mbf u)(t)}. \label{eqn:picard_gradient} \vspace{-3.5mm}
\end{equation}
 
 
Next, we show that $\{S_n\}_{n=1}^{\infty}$ is a Cauchy sequence in $C([0, T]\times X \times \Theta)$ and hence converges to the gradient of the solution map, $\phi$.

\begin{lem}  \label{lem:picard_sensitivity_is_cauchy}
Suppose $\Gamma \subset \R$, $X\subset\R^{n_x}$, $\Theta\subset\R^{n_p}$ are  compact sets and $f \in C^1(\Gamma\times\R^{n_x}\times \Theta)$ is locally Lipschitz. Let $t_0 \in\Gamma$, $a = 2 \sup_{x \in X} \|x\|_2$ and $T < \min\{\frac{1}{K_x}, \frac{a}{2M} \}$, where $K_x$, $M$ are defined as in Thm.~\ref{thm:picard_contraction}. Let $K_{\theta}$ be as defined in Lem.~\ref{lem:picard_is_lipschitz_in_parameters} and  $\nabla_x f, \nabla_\theta f$ be locally Lipschitz continuous. Then, for any  $\mbf u \in C_a[0, T]$, if $\{S_n\}_{n=1}^{\infty}$ is as defined in Eqn.~\eqref{eqn:picard_gradient}, $\{S_n\}_{n=1}^{\infty}$ is a uniformly bounded Cauchy sequence in $C([0, T] \times X\times\Theta)$ with $\|S_n\|_{\infty} \leq \max\big\{\frac{1}{1-TK_{x}}; \frac{TK_{\theta}}{1-TK_{x}}\big\}$.


\end{lem}\vspace{-3.5mm}
\begin{pf}  
First, we show that $S_n \in C([0, T] \times X\times\Theta)$. Define the function $\mbf v_n(t, x, \theta) := ({\mcl P}_{t_0, x, \theta}^n \mbf u)(t)$. Then, since $f \in C^1(\Gamma\times \R^{n_x} \times\Theta)$ and $\mbf u \in C_a[0, T]$, by inductive application of Lem.~\ref{lem:picard_is_continuous} $\mbf v_n  \in C^1([0, T]\times X \times\Theta)$. Thus, $\nabla_{x, \theta} \mbf v_n \in C([0, T]\times X \times\Theta)$.
Next, we show \vspace{-4mm}
\[
\sup_{t \in [0,T]}\|S_n(t,x,\theta)\|_{\infty} \leq  \max\left\{\frac{1}{1-TK_{x}}; \frac{TK_{\theta}}{1-TK_{x}}\right\}.\vspace{-3mm}
\] 
Since $\mcl P_{t_0, x, \theta}^n \mbf u\in C^1([0, T]\times X \times\Theta)$, the Lipschitz bound in Lem.~\ref{lem:picard_is_lipschitz_in_parameters} implies 
$\big\|\nabla_{x,\theta} (\mcl P_{t_0, x, \theta}^n \mbf u)(t)\big\|_2 \leq \max\big\{\frac{1}{1-TK_{x}}; \frac{TK_{\theta}}{1-TK_{x}} \big\}$
for all $(t,x,\theta) \in [0,T] \times X \times \Theta$. 


Finally, we show that $\{S_n\}_{n=1}^{\infty}$ is Cauchy. From Lem.~17 (See Appendix~B) there exist $K_{s}, M_s \geq 0$ such that  for any $n,m\in \N$, $n \geq m$, $(x, \theta) \in X\times \Theta$  and $\mbf u \in C_a[0, T]$  we have \vspace{-3.5mm}
%
%
\begin{align*}
& \|\nabla_{x, \theta} {\mcl P}_{t_0, x, \theta}^n \mbf u - \nabla_{x, \theta} {\mcl P}_{t_0, x, \theta}^m \mbf u\|_\infty  \\
&\quad\leq \hspace{-1mm}(TK_{x})^{m} \big[ mK_{s} \|{\mcl P}_{t_0, x, \theta}^{n-m}\mbf u- \mbf u\|_\infty +  \| \nabla_{x, \theta}  \mcl P^{n-m}_{t_0, x, \theta} \mbf u \|_\infty \big].  \\[-8mm]\notag
\end{align*}  

Furthermore, $\mbf u\in C_a[0, T]$ and Thm.~\ref{thm:picard_contraction} imply $\mcl P^{n-m}_{t_0, x, \theta} \mbf u \in C_a[0, T]$ and there exists $M_s > 0$ such that $\|\nabla_{x, \theta}  \mcl P^{n-m}_{t_0, x, \theta} \mbf u\|_{\infty} \leq M_s$. Hence, we have \vspace{-3.5mm}
%
%
\begin{align*}
&\sup_{t\in[0, T]} \hspace{-2mm}\|\nabla_{x, \theta} (\mcl P_{t_0, x, \theta}^n \mbf u)(t) - \nabla_{x, \theta} (\mcl P_{t_0, x, \theta}^m \mbf u)(t)\|_2  \\[-2mm]
& \qquad\leq (TK_{x})^{m} \big[ mK_{s} 2a + M_s \big].  \\[-8mm]\notag
\end{align*}
Since this bound holds uniformly on $x,\theta\in X \times \Theta$ we have that $\{S_n\}_{n=1}^{\infty}$ is Cauchy.$\blacksquare$  
\end{pf}\vspace{-3.5mm}
 
An immediate consequence of Lem.~\ref{lem:picard_sensitivity_is_cauchy} is that $\{S_n\}_{n=1}^{\infty}$ converges to the gradient of the solution map.  
\begin{prop}\label{prop:picard_sensitivity} 
Suppose that the conditions of Lem.~\ref{lem:picard_sensitivity_is_cauchy} are satisfied and  $\phi$ is the solution map of the ODE defined by $f$.
 Then for all  $\mbf u \in C_a[0, T]$ and $(t,x,\theta) \in [0, T]\times  X\times \Theta$, we have $\lim\limits_{n\rightarrow \infty} \nabla_{x,\theta} (\mcl P_{0,x, \theta}^n \mbf u)(t) = \nabla_{x, \theta} \phi(t, x, \theta)$.\vspace{-3.5mm}

\end{prop}\vspace{-.5mm}
\begin{pf} 
Since for all $\mbf u \in C_a[0, T]$ iterations $S_n$ are a Cauchy sequence, then $S_n(t, x, \theta)$ converges to some $S(t, x, \theta)$ uniformly. By Thm. 7.17 in~\cite{rudin1964principles} we have that for all $t \in [0, T]$, $x \in X$, $\theta\in \Theta$  and for all $\mbf u \in C_a[0, T]$ \vspace{-3.5mm}
 \begin{align*}
  \lim_{n \rightarrow \infty} \hspace{-1.5mm}\nabla_{x, \theta} (\mcl P_{ 0,x, \theta}^n \mbf u)(t) &\hspace{-1mm}=\hspace{-1mm}  \nabla_{x, \theta}\hspace{-1.5mm} \lim_{n \rightarrow \infty}\hspace{-1mm} (\mcl P_{0, x, \theta}^n \mbf u)(t)\hspace{-1mm}=\hspace{-1mm}  \nabla_{x, \theta} \phi(t, x, \theta),  \\[-8mm]\notag
 \end{align*} 
 where $\phi(t, x, \theta)$ is the solution map. $\blacksquare$  
\end{pf}\vspace{-3.5mm} 

Finally, we show that $\nabla_{x, \theta} \mcl P_{t_0, x, \theta}$ is Lipschitz continuous, thereby satisfying Condition 2) of Thm.~\ref{thm:gradient-contractive-map-is-contraction} (wherein through some abuse of notation, $x \in X$ becomes $(x,\theta) \in X \times \Theta$ and $U$ becomes $C_a[0, T]$). 

\begin{lem}\label{lem:picard_sensitivity_is_continuous}

Suppose that the conditions of Lem.~\ref{lem:picard_sensitivity_is_cauchy} are satisfied. 
Then $\nabla_{x, \theta} \mcl P_{t_0, x, \theta} \mbf u$ is  Lipschitz continuous on $(x,\theta) \in X\times \Theta$ and $\mbf u \in C_a[0, T]$ and there exist  $q < 1$, $K > 0$ and $N \in \N$ such that for $n \ge N$, $(x_1, \theta_1), (x_2, \theta_2)\in X \times \Theta$ and $\mbf u_1, \mbf u_2 \in C_a[0, T]$ \vspace{-3.5mm}
\begin{align*}
&\|\nabla_{x, \theta} \mcl P^n_{t_0,  x_1, \theta_1} \mbf u_1  - \nabla_{x, \theta} \mcl P^n_{t_0,  x_2, \theta_2} \mbf u_2\|_\infty \\
 & \qquad \leq q^n \|\mbf u_1 - \mbf u_2\|_{\infty} + K \| x_1 -  x_2\|_{X} + K\|\theta_1 - \theta_2\|_{\Theta}.    \\[-8mm]\notag
\end{align*} 
%
\end{lem}
\begin{pf}
Lem.~17 in  Appendix~\ref{apd:Sensitivity} implies that, there exists $K_s \geq 0$ such that for any $n \in \N$, $(x, \theta) \in X\times \Theta$  and $\mbf u_1, \mbf u_2 \in C_a[0, T]$  we have\vspace{-3.5mm}
\begin{equation*}
\|\nabla_{x, \theta} {\mcl P}_{t_0, x, \theta}^n \mbf u_1 \hspace{-0.5mm} -\hspace{-0.5mm} \nabla_{x, \theta} {\mcl P}_{t_0, x, \theta}^n \mbf u_2\|_\infty\hspace{-1mm} \leq \hspace{-1mm}n (TK_{x})^{n} \hspace{-0.5mm} K_{s} \hspace{-0.5mm}\|\mbf u_1 \hspace{-0.5mm}-\hspace{-0.5mm} \mbf u_2\|_\infty. \vspace{-3.5mm}
\end{equation*} 
Now, since $q/(TK_x)>1$, there exists $N \in \N$ such that $nK_s(T K_x)^n\le q^n$ for all $n \ge N$
and hence\vspace{-3.5mm} 
\begin{equation*}
\|\nabla_{x, \theta} \mcl P_{t_0, x, \theta}^n \mbf u_1 - \nabla_{x, \theta} \mcl P_{t_0, x, \theta}^n \mbf u_2\|_\infty \leq q^n \|\mbf u_1 - \mbf u_2\|_\infty, \vspace{-3.5mm}
\end{equation*}
which holds uniformly on $(x, \theta) \hspace{-1mm}\in \hspace{-1mm} X \hspace{-1mm}\times \hspace{-1mm}\Theta$ and  $\mbf u_1, \mbf u_2 \in C_a[0, T]$. 

Next, from Lem.~18 in  Appendix~B, for all  $(x_1, \theta_1)$, $(x_2, \theta_2) \in X\times \Theta$ and $\mbf u \in C_a[0, T]$ we have \vspace{-3.5mm}
\begin{align*}
&\|\nabla_{x, \theta}  \mcl P^n_{t_0,  x_1, \theta_1} \mbf u - \nabla_{x, \theta}  \mcl P^n_{t_0,  x_2, \theta_2} \mbf u\|_\infty \\[-1mm]
& \qquad\qquad \leq K\|x_1 - x_2\| + K \|\theta_2 - \theta_1\|_2.   \\[-9mm]\notag
\end{align*}   
We conclude that for $n >N$, \vspace{-3.5mm}
\begin{align*}
&\|\nabla_{x, \theta} \mcl P^n_{t_0,  x_1, \theta_1} \mbf u_1  - \nabla_{x, \theta} \mcl P^n_{t_0,  x_2, \theta_2} \mbf u_2\|_\infty \\[-1mm]
 & \quad \leq q^n \|\mbf u_1 - \mbf u_2\|_{\infty} + K \| x_1 -  x_2\|_{2} + K\|\theta_1 - \theta_2\|_2. \blacksquare\\[-8mm]\notag
\end{align*} 
\end{pf}\vspace{-1mm} 
To conclude, we have shown that $\mcl P_{t_0, x, \theta}$ satisfies conditions 1) and 2) of Thm.~\ref{thm:gradient-contractive-map-is-contraction}. Next, in Sec.~\ref{sec:conditions_for_Algs}, we give conditions 
under which $L_{\lambda,n}$ (as defined in Eqn.~\eqref{eqn:extended_loss}) likewise satisfies the conditions of Thm.~\ref{thm:gradient-contractive-map-is-contraction} and use Thm.~\ref{thm:gradient-contractive-map-is-contraction} to show that Algs.~\ref{alg3:modified_picard} and~\ref{alg4:ext_picard} converge.



\section{Conditions for Convergence of Gradient-Contractive Algorithm} \label{sec:conditions_for_Algs}
In this section, we give conditions under which Algs.~\ref{alg3:modified_picard} and~\ref{alg4:ext_picard} converge to a fixed point. 
Specifically, we consider Algs.~\ref{alg3:modified_picard} to be a special case of Alg.~\ref{alg4:ext_picard} and use the results of Sec.~\ref{sec:picard} to propose conditions on the loss function, $L_{\lambda, n}$ (as defined in~\eqref{eqn:extended_loss}), under which the conditions of Thm.~\ref{thm:gradient-contractive-map-is-contraction} are satisfied. 
%

First, recall that the gradient-contractive algorithm for step sizes $\alpha > 0$, $\sigma \in (0, 1]$ and order $n\in \N$ is defined by the sequence $(x_{k},\mbf u_{k})=\mcl T^k (x_{0},\mbf u_{0})$ where we say $(x_{k+1},\mbf u_{k+1})=\mcl T (x_{k},\mbf u_{k})$ if \vspace{-3.5mm}
\begin{align*} 
x_{k+1} &=  \Pi_X\big[x_k - \alpha \nabla_x L(x_k, \mcl P^n(x_k, \mbf u_k))\big] \\
\mbf u_{k+1} &=   (1 - \sigma) \mbf u_{k} + \sigma \mcl P(x_{k+1}, \mbf u_k).  \\[-8mm]\notag
\end{align*}

Note that Alg.~\ref{alg4:ext_picard} is a form of gradient-contractive algorithm, with step sizes $\alpha > 0$, $\sigma \in (0, 1]$ and order $n\in \N$, wherein $x\in X$ becomes $(x,\theta) \in X^{\otimes J} \times \Theta$ and $\mbf u \in U$ becomes $\{\mbf u_j\}_{j=1}^J \in C[0,T]^{\otimes J}$ and 
 where we say $([x_{k+1},\theta_{k+1}],\mbf u_{k+1})=\hat{\mcl T}([x_{k},\theta_k],\mbf u_{k})$ if \vspace{-3.5mm}
%
\begin{align}\label{eqn:picard-gradient-contractive-iterations}
 x_{k+1, j} &= \Pi_{X}[x_{k, j} - \alpha \nabla_{x_{k, j}} L_{\lambda, n}(x_{k}, \theta_{k}, \mbf u_{k})] \notag\\
 \theta_{k+1} &= \Pi_{\Theta}[\theta_{k} - \alpha \nabla_{\theta} L_{\lambda, n}(x_{k}, \theta_{k}, \mbf u_{k})] \\
\mbf u_{k+1, j} &= (1 - \sigma) \mbf u_{k, j} + \sigma \mcl P_{(j-1)T,  x_{k+1, j}, \theta_{k+1}} \mbf u_{k, j}\notag. \\[-8mm]\notag 
\end{align}

We now show that if $f,g$ are smooth and  $L_{\lambda,n}$ are strongly convex, Alg.~\ref{alg4:ext_picard} converges to a fixed point for some $\sigma, \alpha$.

%

\begin{thm}\label{thm:ext_main_theorem} 
 Let $X\subset\R^{n_x}$, $\Theta\subset\R^{n_p}$ be compact and convex and suppose $L_{\lambda, n}(x, \theta, \mbf u)$ (as defined in~\eqref{eqn:extended_loss} for some $J\in \N$ and $\lambda >0$) is strongly convex in $(x,\theta) \in X^{\otimes J}\times  \Theta$, uniformly in $\mbf u$ and $m\in \N$, and that $\nabla_x f, \nabla_\theta f, \nabla_x g, \nabla_\theta g$ are Lipschitz continuous .

Let $\hat{\mcl T}$ be as in Eqn.~\eqref{eqn:picard-gradient-contractive-iterations}. Then,
there exist $\alpha > 0$, $\sigma \in (0, 1]$, $n\in \N$, $\nu < 1$  and $([x^*,\theta^*], \mbf u^*)$ such that for any $(x_0, \theta_0) \in X^{\otimes J}\times \Theta$ and $\mbf u_0 \in C_a[0, T]^{\otimes J}$ (where $T$ and $a$ are defined in Thm.~\ref{thm:picard_contraction}),\vspace{-3.5mm}
\begin{align*}
& \|{x}_k - {x}^*\|_2 + \|\theta_k - \theta^*\|_2 + \|{\mbf u}_k - {\mbf u}^*\|_{\infty}  \\[-1mm]
&\qquad \leq \nu^k(\|{x}_0 - {x}^*\|_2 + \|\theta_0 - \theta^*\|_2 + \|{\mbf u}_0 - {\mbf u}^*\|_{\infty}),  \\[-8mm]\notag
\end{align*}  
where $([x_k,\theta_k],\mbf u_k)=\hat{\mcl T}^k ([x_0,\theta_0],\mbf u_0)$. 
\end{thm}\vspace{-3.5mm}
\begin{pf} 
Recall that $L_{\lambda, n}$ is defined as \vspace{-5mm}
\begin{align*}
& L_{\lambda, n}(x, \theta, \mbf u) \hspace{-1mm}:= \hspace{-1mm} \sum_{j=1}^{J-1} \lambda \|\mcl P^n_{(j-1)T, x_j, \theta} \mbf u_{j}(T) - x_{j+1}\|_2^2\\[-3mm]
&\hspace{-3mm} +\hspace{-1mm}\frac{1}{2N_s}\hspace{-1mm}  \sum_{i=1}^{N_s}\hspace{-1mm}\|y_i \hspace{-0.5mm}-\hspace{-0.5mm} g(t_i, (\mcl P^n_{(j-1)T, x_j, \theta} \mbf u_{\lfloor t_i/T \rfloor}+1)(t_i\hspace{-0mm}\bmod{T}), \theta)\|_2^2 \notag  \\[-8mm] \notag
\end{align*} 

Let us define the extended Picard iterations $\hat{\mcl P}_{x, \theta}:C_a[0, T]^{\otimes J} \rightarrow C_a[0, T]^{\otimes J}$ as $\big(\hat{\mcl P}_{x, \theta} \mbf u\big)_j = \mcl P_{(j-1)T, x_j, \theta} \mbf u_j$ where for $j\in \overline{1,J}$, we denote by $\mbf u_j\in C_a[0, T]$ the $j^{th}$ element of $\mbf u \in C_a[0, T]^{\otimes J}$. 
Then, the map $\hat {\mcl P}$ allows us to represent the mapping $\hat{\mcl T}$ as in Eqn.~\eqref{eqn:contractive_map} -- i.e. $([x_{k+1},\theta_{k+1}],\mbf u_{k+1})=\hat{\mcl T}([x_{k},\theta_k],\mbf u_{k})$ if \vspace{-3.5mm}
\begin{align*}
 x_{k+1, j} &= \Pi_{X}[x_{k, j} - \alpha \nabla_{x_{k, j}} L_{\lambda, n}(x_{k}, \theta_{k}, \mbf u_{k})] \notag\\
 \theta_{k+1} &= \Pi_{\Theta}[\theta_{k} - \alpha \nabla_{\theta} L_{\lambda, n}(x_{k}, \theta_{k}, \mbf u_{k})] \\[-1mm]
\mbf u_{k+1} &= (1 - \sigma) \mbf u_{k} + \sigma \hat{\mcl P}_{x_{k+1}, \theta_{k+1}} \mbf u_k,\\[-8mm]
\end{align*}
where $\mbf u_{k+1} \in C_a[0, T]^{\otimes J}$ and hence $\mbf u_{k+1,j} \in C_a[0, T]$ as specified above.

In the rest of the proof Thm.~\ref{thm:picard_contraction} and Lemmas~\ref{lem:picard_is_lipschitz_in_parameters} and~\ref{lem:picard_sensitivity_is_continuous} are used to show that $\hat{\mcl P}_{x_k, \theta_k}$ and $L_{\lambda, n}$ satisfy the conditions of Thm.~\ref{thm:gradient-contractive-map-is-contraction}.

First, we define  $\hat U:= C_a[0, T]^{\otimes J}$ and $\hat X = X^{\otimes J}$ (with norms $\norm{\hat x}_{\hat X}=\max_i \norm{\hat x_i}_2$ and $\norm{\hat{\mbf u}}_{\hat U}=\max_i \norm{\hat{\mbf u}_i}_\infty$, respectively). Then, since $X$ is compact, $\hat X$ is also compact.  Furthermore, since $a = 2\sup_{x \in X} \|x\|$ we have \vspace{-3.5mm}
\[
2\sup_{\hat x \in \hat X} \norm{\hat x}_{\hat X}=2\max_{j \in \overline{1, J}} \; \; \sup_{\hat x_j \in X} \norm{\hat x_j}_X= a.\vspace{-4mm}
\]


Define  $\hat f \in C([0, T]\times\R^{J \cdot n_x}  \times\Theta)$
 as $\hat f(t, \hat x, \theta) = \bmat{\hat f^T_1(t, \hat x, \theta) &  \cdots & \hat f^T_J(t, \hat x, \theta)}^T$, where $\hat f_j(t, \hat x, \theta) = f(t + (j-1)T, \hat x_j, \theta)$. Then\vspace{-4.5mm}
\[
(\hat{\mcl P}_{\hat x, \theta} \mbf u)(t) = \hat x + \int_0^t \hat f(s, \mbf u(s), \theta) ds\vspace{-3.5mm}
\]
 for all $\hat x \in \hat X$, $\theta \in \Theta$, $\mbf u \in \hat U$ and $t \in [0,T]$.
 
 Furthermore,  \vspace{-3.5mm}
\begin{align*} 
  &\sup_{\substack{\norm{\hat x}\le a , \theta\in \Theta\\ t\in [0, T]}} \hspace{-4mm} \|\hat f( t, \hat x, \theta)\|_{\hat X} = \max_{j \in \overline{1, J}}\;  \sup_{\substack{\norm{\hat x_j}\le a, \theta\in \Theta\\ t\in [0, T] }} \hspace{-4mm}\|\hat f_j( t, \hat x, \theta)\|_2 = M \\[-8mm]
    \end{align*}
    and \vspace{-8mm}
    \begin{align*} 
   &\quad\sup_{\substack{\norm{\hat x}\le a,\norm{\hat y}\le a, \theta\in \Theta \\ t\in[0, T] }} \frac{\|\hat f(t, \hat x, \theta) - \hat f(t, \hat y, \theta)\|_{\hat X}}{\|\hat x - \hat y\|_{\hat X}} \\[-3mm]
 &\qquad\qquad = \hspace{-4mm} \sup_{\substack{\norm{\hat x}\le a,\norm{\hat y}\le a,\\ t\in[0, T], \theta\in\Theta}} \hspace{-4mm} \frac{\max_{j \in \overline{1, J}}\|\hat f_j(t, \hat x_j, \theta) - \hat f_j(t, \hat y_j, \theta)\|_X}{\max_{j \in \overline{1, J}}\|\hat x_j - \hat y_j\|_X} \\[-2mm]
 & \qquad\qquad\leq \sup_{\substack{\norm{\hat x}\le a,\norm{\hat y}\le a,\\ t\in[0, T], \theta\in\Theta}} \hspace{-2mm} \frac{\max_{j \in \overline{1, J}} K_x \|\hat x_j-\hat y_j\|_X}{\max_{j \in \overline{1, J}}\|\hat x_j - \hat y_j\|_X} = K_x.  \\[-8mm]\notag
  \end{align*}
Thus, we conclude that, $\hat f$, $\hat \Gamma:=[0, T]$, $\hat X$, and $\Theta$ satisfy the conditions of Thm.~\ref{thm:picard_contraction} with $t_0 =0$ and $a,T,K_x,M,C_a[0, T]$ as defined therein. Therefore, $\hat{\mcl P}_{x, \theta}$ is a contraction in $\mbf u$ uniformly in $x$.
Furthermore, Lem.~\ref{lem:picard_is_continuous} establishes differentiability of $\hat{\mcl P}_{\hat x, \theta}$. This implies Condition 1) of Thm.~\ref{thm:gradient-contractive-map-is-contraction} is satisfied.

Second, since $f \in C^1([0,JT]\times\R^{n_x}\times \Theta)$ and $\nabla_{x, \theta} f$ are Lipschitz continuous, we have that $\hat f \in C^1([0, T]\times\R^{J \cdot n_x}\times \Theta)$ and $\nabla_{x, \theta} \hat f$ are Lipschitz continuous. Thus, since $\hat f$ satisfies the conditions of Lem.~\ref{lem:picard_sensitivity_is_cauchy}, $\nabla_{x, \theta} \hat{\mcl P}_{x, \theta} {\mbf u}$ is Lipschitz continuous with respect to $x, \theta$. Furthermore, by Lem.~\ref{lem:picard_sensitivity_is_continuous} there exist $q < 1, K > 0$, $N \in \N$ such that for all $n \ge N$,  $( x_1, \theta_1), (x_2, \theta_2)\in \hat X \times \Theta$ and ${\mbf u}_1, {\mbf u}_2 \in \hat U$, we have \vspace{-3.5mm}
\begin{align*}
&\|\nabla_{ x, \theta} \hat{\mcl P}^n_{0, x_1, \theta_1} {\mbf u}_1  - \nabla_{ x, \theta} \hat{\mcl P}^n_{0,   x_2, \theta_2} {\mbf u}_2\|_\infty \\
 & \qquad \leq q^n \|{\mbf u}_1 - {\mbf u}_2\|_{\infty} + K \| x_1 - x_2\|_{\hat X} + K\|\theta_1 - \theta_2\|_2.  \\[-8mm]\notag
\end{align*} 
  Thus, $\nabla_{ x, \theta} \hat{\mcl P}_{x, \theta}{\mbf u}$ satisfies Condition 2) of Thm.~\ref{thm:gradient-contractive-map-is-contraction}.


Finally, for given $\lambda\hspace{-1mm}>\hspace{-1mm}0$, if we define $\hat t_i := t_i\hspace{-0.5mm}\bmod{T}$ and \vspace{-3mm}
\[
L_i(x, \theta, \mbf u)\hspace{-0.5mm} = \hspace{-0.5mm}\begin{cases} \hspace{-0.5mm}\frac{1}{2N_s} \|y_i\hspace{-0.5mm} -\hspace{-0.5mm} g(t_i, \mbf u_{\lfloor t_i/T \rfloor +1}(\hat t_i), \theta)\|^2 &\hspace{-3mm}  \text{ if } i \leq N_s \\
\lambda \|x_{i-N_s + 1}\hspace{-0.5mm} -\hspace{-0.5mm}\mbf u_{i - N_s}(T)\|^2 & \hspace{-13mm}\text{ if } N_s \hspace{-1mm} <\hspace{-1mm}  i\hspace{-1mm}  < \hspace{-1mm}  N_s\hspace{-1mm} + \hspace{-1mm}J.
\end{cases}\vspace{-3mm}
\]
Then, $L_{\lambda,n}(x,\theta,\mbf u) =  \sum_{i=1}^{N_s + J-1}L_i(x, \theta, \hat{\mcl P}^n_{0, x, \theta}\mbf u)$. Furthermore,  since $\nabla_{x, \theta} g$ is Lipschitz continuous, $\nabla_{x,\theta} L_i(x, \theta, \mbf u)$ are Lipschitz continuous, and hence Condition 3) of Thm.~\ref{thm:gradient-contractive-map-is-contraction} is satisfied.

Finally, by assumption, $L_{\lambda,n}$ is strongly convex on $(x,\theta) \in X \times \Theta$, uniformly in ${\mbf u} \in \hat U$ and $n \in \N$, and hence Condition 4) of Thm.~\ref{thm:gradient-contractive-map-is-contraction} is satisfied.



We conclude that convexity of $X, \Theta$ and $C_a[0, T]$ imply convexity of $\hat X\times \Theta$ and $\hat U$, respectively, and hence all conditions of Thm.~\ref{thm:gradient-contractive-map-is-contraction} are satisfied wherein $U \mapsto \hat U$ and $X \mapsto \hat X \times \Theta $ are convex and compact, $\mcl P( x, {\mbf u}) \mapsto \hat{\mcl P}_{0, \hat x, \theta}{\mbf u}$, and $\mcl T(x,\mbf u) \mapsto \hat{\mcl T}([x,\theta],\mbf u)$ -- which concludes the proof.$\blacksquare$ 
 
\end{pf}\vspace{-4mm}
\textbf{Remark 1.} Note that in practice we have four hyper-parameters $\alpha, \sigma, T$ and $n$. One option to choose these parameters is to fix $\sigma = \min\{1, \frac{\mu}{K_L - \mu}\}$ and $\alpha = \frac{2}{\mu + K_L}$, where $\mu$ is a strong convexity modulus of $L_{\lambda, n}$ and $K_L$ is Lipschitz constant of $\nabla_{x, \theta} L_{\lambda, n}$. Then, for a fixed $T$, $n$ may be chosen sufficiently large. Alternately, for a fixed $n \geq 1$ we may choose $T$ sufficiently small.

\begin{figure*}[t]
 \centering
    \begin{subfigure}[t]{0.25\textwidth}
        \centering
\includegraphics[width=0.9\textwidth]{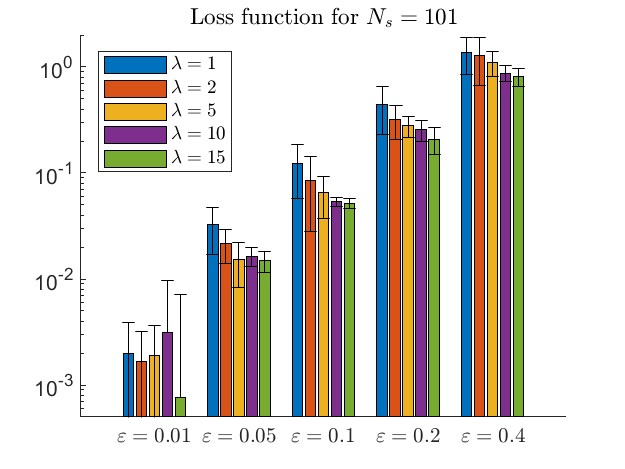}

        \caption{Loss  for $N_s= 101$.}
    \end{subfigure}%
    ~
    \centering
    \begin{subfigure}[t]{0.25\textwidth}
        \centering
\includegraphics[width=0.9\textwidth]{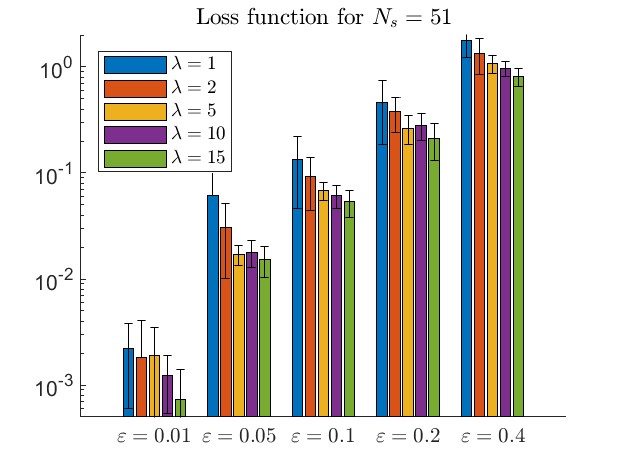}

        \caption{Loss  for $N_s= 51$.}
    \end{subfigure}
    \begin{subfigure}[t]{0.25\textwidth}
        \centering
\includegraphics[width=0.9\textwidth]{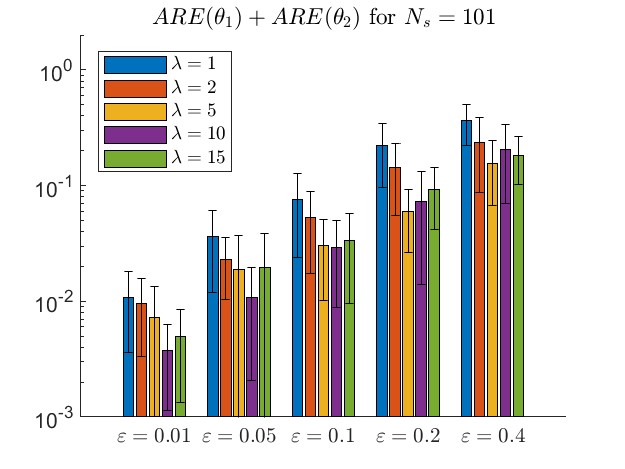}

        \caption{\hspace{-1.0mm}Parameters Errs.,\hspace{-0.5mm}$N_s\hspace{-0.5mm}=\hspace{-0.5mm} 101$.}
    \end{subfigure}%
    ~
    \centering
    \begin{subfigure}[t]{0.25\textwidth}
        \centering
\includegraphics[width=0.9\textwidth]{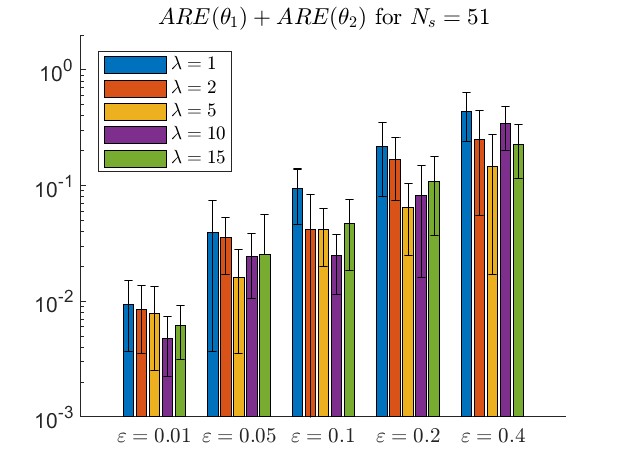}

        \caption{\hspace{-1.0mm}Parameters Errs.,\hspace{-0.5mm}$N_s \hspace{-0.5mm}=\hspace{-0.5mm} 51$.}
    \end{subfigure}\vspace{-2mm}

    \caption{{\footnotesize Loss function as in Eqn.~\eqref{eqn:optimization_problem_unconstrained} and normalized error of parameters ($ARE$) for Van der Pol Oscillator in Eqn.~\eqref{eqn:example_vdp_model} identified by Alg.~\ref{alg4:ext_picard} for $n = 1$, $T = 0.25$ and  $\lambda \in \{1, 2, 5, 10, 15\}$. The true parameters, $\theta_1 = \theta_2 = 1$, variance of the measurements is $\varepsilon \in \{0.01, 0.05, 0.1, 0.2, 0.4\}$. The number of measurements is $N_s = 101$ or $N_s = 51$. The data for $N_s = 51$ is a subsample of the measurements $N_s = 101$ for twice larger sampling times. Average and standard deviation of $ARE$ of parameters and loss function is based on $10$ trials with randomized initial conditions. } }\label{fig:VDP_results}
\end{figure*}

\section{Numerical experiments.} \label{sec:examples}
Having defined a gradient contractive algorithm based on the Picard mapping (Alg.~\ref{alg4:ext_picard}), and having established conditions for convergence of this algorithm to a solution of the parameter estimation problem, we now proceed to evaluate practical aspects of the performance of this algorithm. { This analysis consists of five tests, each focusing on one metric of performance of the algorithm, and each using a nonlinear model and a data set indicative of that metric of performance. For our numerical results, the metric of performance is the average relative error (ARE) in the parameter estimates. Specifically, given $N_{exp}$ estimates of the parameter values, ARE is defined as\vspace{-3.5mm}
 \begin{equation}\label{eqn:ARE}
 ARE = \frac{1}{N_{exp}}.\sum_{i = 1}^{N_{exp}}\frac{| \hat\theta_{i} - \theta|}{|\theta|},\vspace{-3.5mm}
 \end{equation}
where $\hat \theta_i$ are the estimated parameters for the $i^{th}$ data set and $\theta$ is the true parameter values.}

In all examples, initial conditions are unknown, and measurements include varying levels of noise. All tests of Alg.~\ref{alg4:ext_picard} use step size $\sigma=.5$ and a backtracking line search to determine step size $\alpha$. The parameter $T$ varies by example and termination occurs when the norm of gradient $\norm{\nabla_{\theta}L_{\lambda,n}}$ is sufficiently small.  The numerical tests are summarized as follows.

 First, in Subsec.~\ref{subsec:VDP_reg}, we evaluate the effect of the regularization parameter $\lambda$ on error in the parameter estimates using various levels of noise and as applied to the Van der Pol Oscillator, with direct measurement of both states.       
Second, in Subsec.~\ref{subsec:example_FHN_irreg_sampling}, we evaluate the effect of irregular sampling intervals on error in parameter estimates and as illustrated by the  FitzHugh-Nagumo neuron using measurements of a single state.  
Third, in Subsec.~\ref{subsec:example_preypredator_our_vs_others}, we consider performance when data is obtained using multiple instances but with sparse data (long periods between samples) as illustrated by the Rosenzweig-MacArthur model.   
 Fourth, in Subsec.~\ref{subsec:example_tumor_inputs} we consider experiments driven by external excitation and with a large number of states and parameters. We evaluate the effect of different excitations on estimation accuracy as illustrated by a model of tumor growth. 
Finally in Subsec.~\ref{subsec:example_lorenz_multi_trajectories}, we compare the accuracy of the algorithm with both comparable heuristics~\cite{brunton2016discovering} and gradient-free optimization~\cite{sivanandam2008genetic,ljung2012version} as illustrated using the Lorenz model. { Additional numerical examples for estimating unknown parameters in real-world models of the attitude dynamics can be found in~\cite{talitckii2026picard}.}
\subsection{The Van der Pol Oscillator and the Effect of Regularization Parameter}\label{subsec:VDP_reg}

In the first numerical example, we investigate the performance of Alg.~\ref{alg4:ext_picard} (defined by regularized loss function $L_{\lambda, n}$ in Eqn.~\eqref{eqn:extended_loss}) as a function of regularization parameter $\lambda$. Recall, this regularization parameter $\lambda$ weights a penalty for discontinuity of the approximated solution at times $jT$. Hence, for a larger weight $\lambda$, the estimated solution, $\mbf u$ is closer to the solution map on extended time intervals and hence better captures the relationship between estimated parameter values and evolution of the state. However, as $\lambda$ increases, convergence time increases and less weight is placed on matching the estimated solution to the given data. To provide guidance on choice of regularization parameter, we consider Van der Pol oscillator. \vspace{-3.5mm}
\begin{align}
    \dot x_1(t) &= \theta_1 x_2(t) \label{eqn:example_vdp_model} \\
    \dot x_2(t) &= -\theta_1 x_1(t) + \theta_2 (1 - x_1^2(t)) x_2(t), \notag  \\[-8mm]\notag
\end{align}
where both states are observable and our goal is to estimate the parameters $\theta_1, \theta_2$. The data sets used in this example are generated from a single instance (with single initial condition $x_0$) over a time interval $[0,t_f]$ using $t_f=10$ and true parameter values $\theta_1 = \theta_2 = 1$, $N_s$ evenly distributed sampling times $t \in \{0:t_f/(N_s-1) i: t_f\}$ and noisy data of the form  $y_i=(1 + n_i)\phi(t_i, x_0, \theta) + m_i$, 
 where $\phi(t_i,  x, \theta)$ are the actual states and $n_i, m_i$ are normally distributed with zero mean and variance $\varepsilon$.

To apply Alg.~\ref{alg4:ext_picard}, we use the measurement data to approximate a bound on the Lipschitz constant as $K_x\sim 4$ which yields a convergence interval of $T = {1/K_x}= .25s$. This yields $J={t_f/T}= 40$ sub-intervals. A single Picard iteration $n=1$ is used.

To evaluate the effect of regularization, we now construct 100 data sets, $Y_{1,\varepsilon,x_0}$ (using $N_s=101$) and $Y_{2,\varepsilon,x_0}$ (using $N_s=51$) for variances $\varepsilon \in \{0.01, 0.05, 0.1, 0.2, 0.4\}$ and 10 randomly selected initial conditions $x_0\in X_{\lambda,\varepsilon}$ (such that $x_0$ is always inside the limit cycle). Alg.~\ref{alg4:ext_picard} is then applied to each data set using regularization parameters $\lambda \in \{1,2,5,10,15\}$ -- yielding parameter estimates $\hat \theta_{\lambda, \varepsilon,x_0}$ and initial condition estimates $\hat x_{\lambda, \varepsilon,x_0}$. 

 

For each choice of $\lambda,\varepsilon$, performance is evaluated by first performing numerical simulation with estimated parameters, $\hat \theta_{\lambda, \varepsilon,x_0}$, and initial conditions, $\hat x_{\lambda, \varepsilon,x_0}$, and then computing the average loss so that if $\phi$ is the solution map, we have \vspace{-3.5mm}
 \[
L(\lambda, \varepsilon) = \frac{1}{10}\hspace{-2mm}\sum_{x_0\in X_{\lambda,\varepsilon}}\hspace{-3mm}\frac{1}{N_s} \sum\nolimits_{i=1}^{N_s}\|\phi(t_i,\hat x_{\lambda, \varepsilon,x_0},\hat \theta_{\lambda, \varepsilon,x_0}) - y_i\|_2^2.\vspace{-3.5mm}
 \]
In addition, { for a parameter, $\theta$, we compute the average relative error ($ARE$) in the parameter estimates as in Eqn.~\eqref{eqn:ARE}, where the summation is given by parameter estimates for different initial conditions $x_0\in X_{\lambda,\varepsilon}$.}

As seen in Fig.~\ref{fig:VDP_results}, increasing the regularization parameter $\lambda$ significantly improves both average loss and ARE for large noise variance ($\varepsilon \in \{0.1, 0.2, 0.4\}$), but does not significantly improve performance with low noise ($\varepsilon \in \{0.01, 0.05\}$). In addition, as expected, more data ($N_s = 101$ vs $N_s = 51$) improved performance in almost every case -- See Fig.~\ref{fig:VDP_results} a, c) and b, d). 

 {  Although Alg.~\ref{alg4:ext_picard} accurately estimates parameters of the model in~\eqref{eqn:example_vdp_model}, it is also possible to have locally identifiable or unidentifiable parameters.  To illustrate the locally identifiable parameters, let's consider a different parameterization of the model in~\eqref{eqn:example_vdp_model}. Specifically, using a variable substitution $\theta_1 = (\rho - 1)^2$, we have that $\rho = 2$ and $\rho = 0$ gives us $\theta_1 = 1$, both choices leading to identical system response in~\eqref{eqn:example_preypredator_model}. Numerically, using Alg.~\ref{alg4:ext_picard} to estimate $\rho$, we observe that initial parameter estimate $\rho_{0}=-1$ results in convergence to $\rho = 0$  and $\rho_{0}=3$ results in convergence to $\rho = 2$. This illustrates the importance of the choice of parameterization of the model to avoid locally identifiable or non-identifiable models~\cite{barreiro2023origins}.

}
 \begin{figure}
    \centering
\begin{subfigure}[t]{0.23\textwidth}
        \centering
\includegraphics[width=1.0\textwidth]{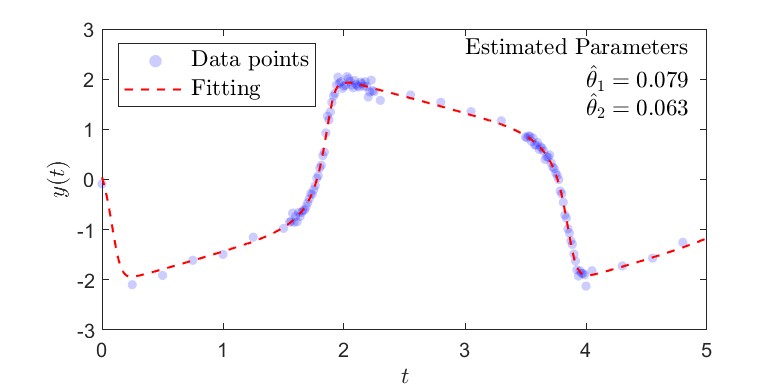}\vspace{-2mm}

        \caption{\hspace{-1.0mm}Results for variance $\varepsilon\hspace{-1mm}=\hspace{-1mm}.05$.}
    \end{subfigure}%
    ~
    \centering
    \begin{subfigure}[t]{0.23\textwidth}
        \centering
\includegraphics[width=1.0\textwidth]{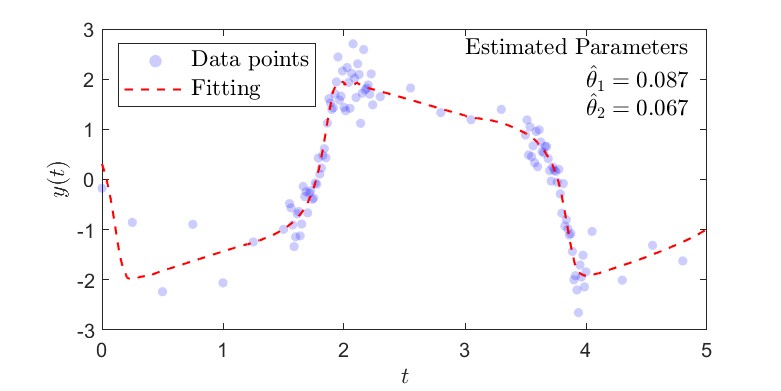}\vspace{-2mm}

        \caption{\hspace{-1.0mm}Results for variance $\varepsilon\hspace{-1mm}=\hspace{-1mm}.2$.}
    \end{subfigure}
 \vspace{-2mm}
    \caption{{\footnotesize Parameter Estimation using Alg.~\ref{alg4:ext_picard} for the FitzHugh-Nagumo model in Eqn.~\eqref{eqn:FHN_example_model} with irregular sampling times (blue dots) clustered around spiking and bursting with noise variances of $\varepsilon=.05$ (Fig 2a) and $\varepsilon=.2$ (Fig 2b). The red dashed line indicates simulation based on the estimated parameters ($\hat \theta_i$) and initial states. True parameter values are $\theta_1 =0.08$ and $\theta_2 = 0.064$. 
    }   }\vspace{-1.5mm}
    \label{fig:FHN_results}
\end{figure}
  
\subsection{The FitzHugh-Nagumo Neuron and the Effect of Irregular Sampling} 
 \label{subsec:example_FHN_irreg_sampling}
\begin{table*}[t]\centering
\begin{tabular}{c|c|cccccc|c}
$\#$ of Initial  & Total Number of & \multicolumn{6}{|c|}{Identified Parameters} & \\
 Conditions & Samples & $ARE(a)$ & $ARE(b)$ & $ARE(c)$ & $ARE(r)$ & $ARE(k)$ & $ARE(K)$ & $ \sum\limits_{i = 1}^{6} \hspace{-1.0mm}\frac{ARE(\theta_i)}{6}$ \\
\hline
 $N_c = 1$ & $N_s N_c = 5$ & $0.19 \hspace{-0.5mm}\pm \hspace{-0.5mm} 0.11$& $0.37 \hspace{-0.5mm} \pm \hspace{-0.5mm} 0.29$& $0.14 \hspace{-0.5mm} \pm \hspace{-0.5mm} 0.12$& $0.35  \hspace{-0.5mm}\pm  \hspace{-0.5mm}0.36$& $0.30  \hspace{-0.5mm}\pm \hspace{-0.5mm} 0.31$& $0.21 \hspace{-0.5mm} \pm  \hspace{-0.5mm}0.15$& \hspace{-1mm} $0.26 \hspace{-0.5mm} \pm \hspace{-0.5mm} 0.18$ \\ 
 $N_c = 2$ & $N_s N_c = 10$ & $0.18 \hspace{-0.5mm} \pm  \hspace{-0.5mm}0.19$& $0.17 \hspace{-0.5mm} \pm \hspace{-0.5mm} 0.11$& $0.10 \hspace{-0.5mm} \pm \hspace{-0.5mm} 0.09$& $0.13 \hspace{-0.5mm} \pm  \hspace{-0.5mm}0.09$& $0.17 \hspace{-0.5mm} \pm  \hspace{-0.5mm}0.15$& $0.07 \hspace{-0.5mm} \pm \hspace{-0.5mm} 0.06$& \hspace{-1mm} $0.14 \hspace{-0.5mm} \pm \hspace{-0.5mm} 0.09$ \\ 
 $N_c = 3$ & $N_s N_c = 15$ & $0.17 \hspace{-0.5mm} \pm  \hspace{-0.5mm}0.15$& $0.17 \hspace{-0.5mm} \pm  \hspace{-0.5mm}0.15$& $0.09  \hspace{-0.5mm}\pm \hspace{-0.5mm} 0.06$& $0.08  \hspace{-0.5mm}\pm \hspace{-0.5mm} 0.09$& $0.14 \hspace{-0.5mm} \pm \hspace{-0.5mm} 0.12$& $0.06 \hspace{-0.5mm} \pm \hspace{-0.5mm} 0.05$& \hspace{-1mm} $0.12 \hspace{-0.5mm} \pm  \hspace{-0.5mm}0.08$ \\ 
 $N_c = 4$ & $N_s N_c = 20$ & $0.13 \hspace{-0.5mm} \pm \hspace{-0.5mm} 0.11$& $0.15 \hspace{-0.5mm} \pm \hspace{-0.5mm} 0.16$& $0.07 \hspace{-0.5mm} \pm \hspace{-0.5mm} 0.06$& $0.06 \hspace{-0.5mm} \pm \hspace{-0.5mm} 0.05$& $0.11 \hspace{-0.5mm} \pm \hspace{-0.5mm} 0.10$& $0.06 \hspace{-0.5mm} \pm \hspace{-0.5mm} 0.05$& \hspace{-1mm} $0.10 \hspace{-0.5mm} \pm \hspace{-0.5mm} 0.07$ \\

\end{tabular}
\caption{\vspace{-3.5mm}{\footnotesize Normalized error in estimated parameters using Alg.~\ref{alg4:ext_picard} for the Rosenzweig-MacArthur predator-prey model~\eqref{eqn:example_preypredator_model} vs. number of populations sampled using sparse sampling times $t\in \{0, 1.25, 2.5, 3.75, 5\}$. Each test is repeated 10 times using randomized initial conditions noise variance $\varepsilon=.05$ and Average Relative Error (ARE) and standard deviation are listed.  True parameter values are $a = 2.8$, $b=0.7$, $c=1.35$, $r = 3.5$, $k = 1.5$, $K = 1.4$. }
}
\label{tab:pp_results}
\end{table*}
In this subsection, we evaluate the effect of irregular sampling intervals on error in estimated parameters. For this case, we consider the FitzHugh-Nagumo model a single cortical neuron, which represents spiking and bursting behavior of neurons~\cite{fitzhugh1961impulses}. This model is a standard example of a stiff ODE and hence it is known that sampling at times of spiking and bursting are more useful than during interim periods -- See Fig.~\ref{fig:FHN_results}. Specifically, we consider the model \vspace{-3.5mm}
\begin{align}
\dot x_1(t) &=  10(x_1(t) - 1/3 x_1^3(t) - x_2(t) + I(t)) \notag \\
\dot x_2(t) &=  10(\theta_1 x_1(t) - \theta_2 x_2(t) + 0.056) \label{eqn:FHN_example_model} \\
y(t) & = g(x(t)):=x_1(t), \notag \\[-8mm]\notag
\end{align}
where the states, $x_1(t), x_2(t) \in \R$, represent membrane potential of neurons and membrane recovery variables and the input, $I(t)$, represents the stimulus current and is set to $I(t) = 1$ when the neuron is active. The measured output is $y(t)$ and the parameters to be estimated are $\theta_1, \theta_2 \geq 0$.

For this test, we use true parameter values $\theta_1 = 0.08$ and $\theta_2 = 0.064$. Based on observations of the output, $y(t)=x_1(t)$, sampling times are clustered around observed spiking and bursting times, $T_2:=\{1.55+0.0125k\}_{k=0}^{56}$ and $T_4:= \{3.5+0.0125k\}_{k=0}^{40}$. Sampling times during recovery periods are less frequent: $T_1:= \{0.25k\}_{k=0}^6$, $  T_3 := \{2.3+0.25k\}_{k=0}^4$, and $T_5:= \{4.05 + 0.25k\}_{k=0}^3$ where the union of all sampling times is defined as $t_i \in T := \cup_{i=1:5} T_i$. These sampling times are depicted as blue dots in Fig.~\ref{fig:FHN_results} and where (as in Subsec.~\ref{subsec:VDP_reg}) measurements are taken over a single instance using initial condition $x_1(0)=0, x_2(0)=2$ as \vspace{-3.5mm}
\[
y_i=(1 + n_i)g(\phi(t_i, x(0), \theta)) + m_i\qquad  t_i \in T,\vspace{-3.5mm}
\]
where $\phi$ is the solution map and $n_i, m_i$ are normally distributed with zero mean and variances of $\varepsilon = .05$ (Fig~\ref{fig:FHN_results}a) and $\varepsilon = .2$ (Fig~\ref{fig:FHN_results}b). 
Applying Alg.~\ref{alg4:ext_picard} with $\lambda =5$, $n =2$, and $T=.1$ we obtain parameter estimates $\hat \theta_1 = 0.079$ and $\hat \theta_2 = 0.063$ for variance $\varepsilon = .05$; and $\hat \theta_1 = 0.087$ and $\hat \theta_2 = 0.067$ for variance $\varepsilon=.2$. 

%
%

\subsection{Rosenzweig-MacArthur predator-prey model and the Effect of Sparse Data.}
 \label{subsec:example_preypredator_our_vs_others}

In the next example, we evaluate performance of Alg.~\ref{alg4:ext_picard} for sparse data using the Rosenzweig-MacArthur predator-prey model~\cite{grunert2021evolutionarily}, wherein the states represent populations of predators and times are measured in years. As in many biological models, collecting data entails the labor-intensive counting of populations and hence data for such models is typically limited to a few samples. In addition, measurements of biological processes are often inaccurate and hence repeated experiments are used to obtain statistically significant sample sizes. For example, in the case of predator-prey models, data may be based on tracking several distinct populations in order to reduce the influence of unmodelled and localized factors such as geography. In our case, the model takes the following form~\cite{grunert2021evolutionarily} \vspace{-3.5mm} 
\begin{align}
 \hspace{0mm}\dot x_1(t) &= rx_1(t)\big(1 - \frac{x_1(t)}{K}\big) -\frac{ ax_1(t)x_2(t)}{b + x_1(t)}\notag  \\[-2mm] 
 \hspace{0mm} \dot x_2(t) &= k \frac{ ax_1(t)x_2(t)}{b + x_1(t)} - c x_2(t), \label{eqn:example_preypredator_model}\\[-9mm]\notag
\end{align}
where $x_1(t), x_2(t)$ represent prey and predator populations, respectively. This model has $6$ unknown parameters ($a,b,c,r,k,K$), with: $r, K$ representing reproduction rate and carrying capacity of the prey population; $a, k$ are encounter and growth rates, $b$ is an environmental factor, and $c$ is predator death rate. For this test, we use true parameters $a = 2.8$, $b=0.7$, $c=1.35$, $r = 3.5$, $k = 1.5$ and $K = 1.4$.

To evaluate the effect of sparse sampling and multiple instances, the algorithm is tested using data obtained from between 1 and 4 distinct populations ($N_c=1:4$ is number of initial conditions), with population counts being made at years $t_i \in \{0, 1.25, 2.5, 3.75, 5\}$ using measurement with noise model $y_{ij} =  (1 + n_i) \phi(t_i, x_j, \theta) + m_i$,
where $n_i, m_i$ are normally distributed with zero mean and variance $\varepsilon = 0.05$.

 Alg.~\ref{alg4:ext_picard} is then applied to each case with parameters $\lambda = 5$, $n = 2$ and $T = 0.2$. The test is repeated ten times and the average relative error in the parameter estimates is listed in Table~\ref{tab:pp_results}, where even for small amount of samples we were able to identify parameters. { As expected, these results demonstrate that increasing the number of sampled populations reduces error in the parameter estimates, with reasonable accuracy achieved using as few as two populations. However, certain parameters ($a$ and $b$) do not appear to converge, even with an increased number of populations. Furthermore, we find that for noise variance $\varepsilon = 0.05$, the objective value for, e.g., $a=3.65$ and $b=0.83$ is indistinguishable from the objective value obtained using the true values of $a=2.8$ and $b=.7$. This implies that the problem here is not the existence of local minima, but rather that the model parameterization may not be identifiable as stated.

 } 
 
%
%
%

\subsection{Tumor Growth, Identifiability, and Excitation.}

 \label{subsec:example_tumor_inputs}
 
In this example, we examine the significance of model excitation for a weakly identifiable model. As discussed in~\cite{baaijens2016existence, raue2014comparison}, biological models with multiple states and parameters are often weakly identifiable in that multiple choices of parameters yield the same or approximately the same solution. We examine here the effect of using datasets generated from multiple non-redundant excitations to improve identifiability of these parameters. Specifically, we use a model of tumour immunodynamics with chemotherapy~\cite{de2003dynamics} where different drug dosing strategies are used as inputs to the model, given as  \vspace{-3.5mm}
\begin{align}
\dot N(t) &= \theta_1 N(t)(1 \hspace{-0.5mm}- \hspace{-0.5mm}  N(t))  \hspace{-0.5mm}- \hspace{-0.5mm} \theta_1 N(t) T(t)  \hspace{-0.5mm}- \hspace{-0.5mm}  (1 \hspace{-0.5mm} - \hspace{-0.5mm} e^{-v(t)}) N(t)  \notag\\
\dot T(t) &= \theta_2 T(t)(1  -  T(t)) - \theta_3 T(t) I(t) - \theta_1 T(t) N(t) \notag\\[-2mm]
& \qquad\qquad\qquad\qquad\qquad\qquad- 3(1 - e^{-v(t)}) T(t)\notag\\[-2mm]
\dot I(t) &= 	\theta_4 + \theta_5 \frac{I(t) T(t)}{\rho + T(t)} - \theta_{1} I(t) T(t) - \theta_{6} I(t) \label{eqn:example_tumor_model}\\[-4mm]
& \qquad\qquad\qquad\qquad\qquad\qquad - 2 (1 - e^{-v(t)}) I(t) \notag \\[-2mm]
\dot v(t) &= -  v(t) + u(t) \notag \\[-2mm]
y(t) &= g(N(t), T(t), I(t), v(t)) := \bmat{N(t) & T(t)  & I(t)}^T, \notag\\[-8mm]\notag
\end{align}
where $N(t)$ is the density of healthy cells, $T(t)$ is tumor cells, $I(t)$ is immune cells and $v(t)$ is the concentration of drugs in the tumor area. The system has a single input $u(t)$, that represents drug dosing. We assume that states $N(t), T(t), I(t)$ are outputs of the system. Parameter values are taken from~\cite{de2003dynamics}, where $\theta_{1} = 1$, $\theta_2 = 1.5$, $\theta_3 = 0.5$, $\theta_4 = 0.33$, $\theta_5 = 0.01$, $\theta_{6} = 0.2$ and $\rho = 0.3$. The initial conditions are $N(0) = 1$, $T(0) = 2$, $I(0) = 1.65$, $v(0) = 0$. 

5 drug dosing strategies are chosen as $u_1(t) = 0$, \vspace{-3.5mm}
\begin{align*} 
 u_2(t) & = \begin{cases} 
  e^{-t}, & \hspace{-0.5mm}\text{if } t \in [0, 3] \\
  0, & \hspace{-0.5mm}\text{otherwise}
 \end{cases} \quad \hspace{-0.5mm}
 u_3(t)   = \begin{cases} 
  e^{t}, &\hspace{-0.5mm}\hspace{-1mm} \text{if } t \in [0, 3] \\
  0, &\hspace{-0.5mm}\hspace{-1mm} \text{otherwise}
 \end{cases}\\
u_4(t) & = \begin{cases} 
  3\hspace{-0.5mm} -\hspace{-0.5mm} t & \hspace{-0.5mm}\text{if } t \in [3, 6] \\
  0, & \hspace{-0.5mm}\text{otherwise}
 \end{cases} \quad \hspace{-0.5mm}\hspace{-0.5mm}\hspace{-0.5mm}
  u_5(t)   = \begin{cases} 
  e^t, &\hspace{-0.5mm}\hspace{-1mm} \text{if } t \in [1.5, 4.5] \\
  0, &\hspace{-0.5mm}\hspace{-1mm} \text{otherwise,}
 \end{cases} \\[-9mm]\notag
\end{align*}
where it is observed that $u_3,u_4,u_5$ clear the tumour, but $u_1$ and $u_2$ do not.

\begin{table}\centering

\begin{tabular}{c|ccccc}
 \hspace{-4mm}{\tiny Identified} \hspace{-3mm} &  \multicolumn{5}{|c}{Available Inputs} \\ 
 \hspace{-4mm}{\tiny Parameters} \hspace{-3mm} & $u_1$ & $u_{1:2}$ & $u_{1:3}$ & $u_{1:4}$ & $u_{1:5}$\\  
\hline
\hspace{-4mm}\tiny$ARE(\theta_1)$\hspace{-3mm} & {\hspace{-1mm}\color{red}\tiny$0.34\pm 0.06$}\hspace{-1mm} & {\hspace{-1mm}\color{UBCgreen}\tiny$0.18\pm 0.05$}\hspace{-1mm} & {\hspace{-1mm}\color{UBCgreen}\tiny$0.23\pm 0.15$}\hspace{-1mm} & {\hspace{-1mm}\color{UBCgreen}\tiny$0.12\pm 0.02$}\hspace{-1mm} & {\hspace{-1mm}\color{UBCgreen}\tiny$0.03\pm 0.03$}\hspace{-1mm} \\  
  
\hspace{-4mm}\tiny$ARE(\theta_2)$\hspace{-3mm} & {\hspace{-1mm}\color{red}\tiny$0.44\pm 0.09$}\hspace{-1mm} & {\hspace{-1mm}\color{UBCgreen}\tiny$0.21\pm 0.07$}\hspace{-1mm} & {\hspace{-1mm}\color{red}\tiny$0.43\pm 0.13$}\hspace{-1mm} & {\hspace{-1mm}\color{UBCgreen}\tiny$0.28\pm 0.04$}\hspace{-1mm} & {\hspace{-1mm}\color{UBCgreen}\tiny$0.07\pm 0.07$}\hspace{-1mm} \\  
  
\hspace{-4mm}\tiny$ARE(\theta_3)$\hspace{-3mm} & {\hspace{-1mm}\color{red}\tiny$0.67\pm 0.16$}\hspace{-1mm} & {\hspace{-1mm}\color{red}\tiny$0.31\pm 0.13$}\hspace{-1mm} & {\hspace{-1mm}\color{red}\tiny$0.88\pm 0.05$}\hspace{-1mm} & {\hspace{-1mm}\color{red}\tiny$0.67\pm 0.08$}\hspace{-1mm} & {\hspace{-1mm}\color{UBCgreen}\tiny$0.15\pm 0.16$}\hspace{-1mm} \\  
  
\hspace{-4mm}\tiny$ARE(\theta_4)$\hspace{-3mm} & {\hspace{-1mm}\color{red}\tiny$0.43\pm 0.09$}\hspace{-1mm} & {\hspace{-1mm}\color{UBCgreen}\tiny$0.20\pm 0.07$}\hspace{-1mm} & {\hspace{-1mm}\color{red}\tiny$0.33\pm 0.18$}\hspace{-1mm} & {\hspace{-1mm}\color{UBCgreen}\tiny$0.17\pm 0.04$}\hspace{-1mm} & {\hspace{-1mm}\color{UBCgreen}\tiny$0.05\pm 0.03$}\hspace{-1mm} \\  
  
\hspace{-4mm}\tiny$ARE(\theta_5)$\hspace{-3mm} & {\hspace{-1mm}\color{red}\tiny$0.78\pm 0.32$}\hspace{-1mm} & {\hspace{-1mm}\color{red}\tiny$2.70\pm 1.50$}\hspace{-1mm} & {\hspace{-1mm}\color{red}\tiny$3.96\pm 0.00$}\hspace{-1mm} & {\hspace{-1mm}\color{red}\tiny$3.96\pm 0.00$}\hspace{-1mm} & {\hspace{-1mm}\color{red}\tiny$1.39\pm 0.99$}\hspace{-1mm} \\  
  
\hspace{-4mm}\tiny$ARE(\theta_6)$\hspace{-3mm} & {\hspace{-1mm}\color{red}\tiny$0.70\pm 0.18$}\hspace{-1mm} & {\hspace{-1mm}\color{red}\tiny$0.37\pm 0.03$}\hspace{-1mm} & {\hspace{-1mm}\color{red}\tiny$0.78\pm 0.18$}\hspace{-1mm} & {\hspace{-1mm}\color{red}\tiny$0.57\pm 0.08$}\hspace{-1mm} & {\hspace{-1mm}\color{UBCgreen}\tiny$0.01\pm 0.01$}\hspace{-1mm} \\  
  
\hspace{-4mm}\tiny$ARE(\theta_7)$\hspace{-3mm} & {\hspace{-1mm}\color{red}\tiny$0.66\pm 0.00$}\hspace{-1mm} & {\hspace{-1mm}\color{red}\tiny$0.69\pm 0.02$}\hspace{-1mm} & {\hspace{-1mm}\color{red}\tiny$0.99\pm 0.00$}\hspace{-1mm} & {\hspace{-1mm}\color{red}\tiny$0.99\pm 0.00$}\hspace{-1mm} & {\hspace{-1mm}\color{red}\tiny$0.96\pm 0.05$}\hspace{-1mm}

\end{tabular}
\caption{ \vspace{-3.5mm}{\footnotesize Normalized error in identified parameters ($\theta_1, ..., \theta_6, \rho$) for model of tumor growth as in~\eqref{eqn:example_tumor_model} based on multiple input-output measurements. Each pair of input-output measurements has $N_s = 51$ data samples and is generated by input $u$ from the set of observed inputs $u_{1:5}$ with noise variance $0.01$. True parameters are $\theta_{1} = 1$, $\theta_2 = 1.5$, $\theta_3 = 0.5$. $\theta_4 = 0.33$, $\theta_5 = 0.01$, $\theta_{6} = 0.2$ and $\rho = 0.3$. The initial condition $N(0)= 1$, $T(0) = 2$, $I(0) = 1.65$, $v(0) = 0$ is identical for each input-output pair. Average and standard deviation of parameter errors are presented based on $3$ different trials. The green values represent that the parameter has been identified ($ARE < 0.3$) and the red values represent that the parameter has not been identified ($ARE > 0.3$).}
}
\label{tab:tumor_results}
\end{table}
For each dosing strategy, $i=1,\cdots,5$, we generate 3 datasets $\{Y_{i,j}\}_{j=1}^3$ each using measurement noises $n_i, m_i$ (normally distributed with zero mean and variance $\varepsilon=.01$) so that $y_{i} = (1 + n_i)g(\phi_u(t_i, x, \theta)) + m_i$ where $\phi_u$ is the solution of Eqn.~\eqref{eqn:example_tumor_model} for input $u(t)$, and $t_i = \{0:\frac{t_f}{N_s - 1}i : t_f\}$, where $N_s = 51$ and $t_f = 6$.  

For testing, we apply Alg.~\ref{alg4:ext_picard} (with $\lambda = 2$, $n = 2$ and $T = 0.1$) to the datasets $Y_{1,j}$, $Y_{1,j}\cup Y_{2,j}$, $Y_{1,j}\cup Y_{2,j}\cup Y_{3,j}$, $Y_{1,j}\cup Y_{2,j}\cup Y_{3,j}\cup Y_{4,j}$, and $Y_{1,j}\cup Y_{2,j}\cup Y_{3,j}\cup Y_{3,j}\cup Y_{5,j}$  -- representing an increase in the number of excitations used to identify the data. In each case we determine the ARE (over $j=1:3$) for each identified parameter. The results are listed in Table~\ref{tab:tumor_results}, where green indicates that the ARE for the given parameter is less that $0.3$. These results indicate the failure to include any excitation ($Y_{1j}$) results in failure to effectively identify the parameters while inclusion of inputs $u_2,u_3,u_4$ generates improved results. {  Inclusion of all 5 excitations generates the best results, although even in this case, it seems that some parameters are not identifiable. This result suggests that the algorithm converges to local optima of the loss function, a typical behavior for parameter estimation problems. Indeed, the achieved objective value for the estimated parameters is $L_{ls} = 0.027$ compared to $L_{ls} = 0.006$ for the true parameter values. Thus, in contrast to Example~\ref{subsec:example_preypredator_our_vs_others}, the proposed algorithm was not able to find a globally-optimal solution. The failure of the convergence may be caused by the non-convexity of the optimization problem or the convergence to approximated KKT conditions.  Thus, to more accurately identify the parameters, reparameterization of the model may be necessary~(Sec.~5 in~\cite{MESHKAT201119} or Sec.~2 in \cite{baaijens2016existence}).}



\begin{figure*}[t!]
 \centering
    \begin{subfigure}[t]{0.32\textwidth}
        \centering
\includegraphics[width=.8\textwidth]{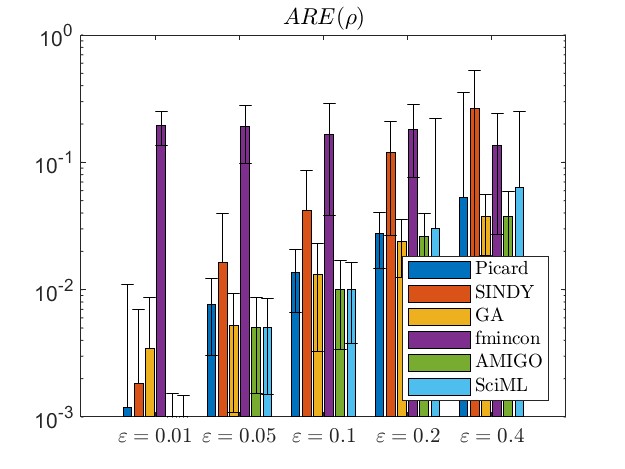}
        \caption{Error in identified $\rho$ parameter.}
    \end{subfigure}%
    ~
    \centering
    \begin{subfigure}[t]{0.32\textwidth}
        \centering
\includegraphics[width=.8\textwidth]{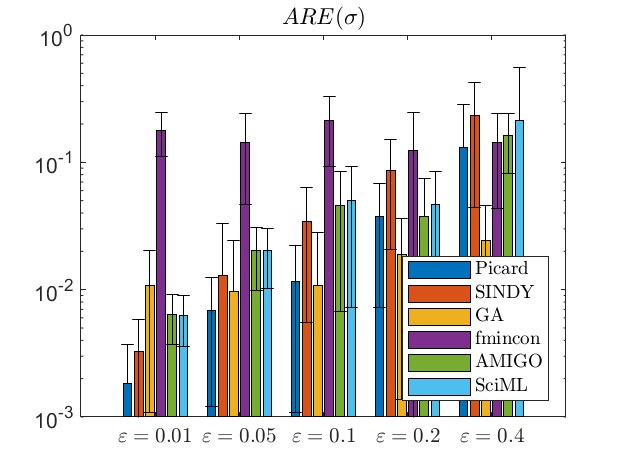}
        \caption{Error in identified $\sigma$ parameter.}
    \end{subfigure}
    \begin{subfigure}[t]{0.32\textwidth}
        \centering
\includegraphics[width=.8\textwidth]{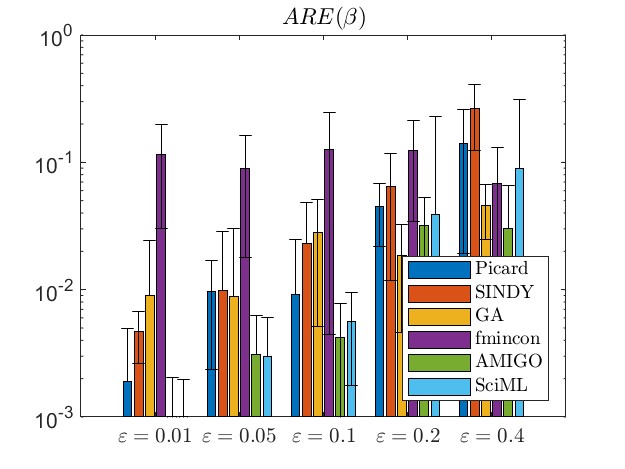}
        \caption{Error in identified $\beta$ parameter.}
    \end{subfigure}%
    \vspace{-2mm}

    \caption{ {\footnotesize   Accuracy in identified system parameters for the Lorenz system (Eqn.~\eqref{eqn:example_lorenz_model}) using Alg.~\ref{alg4:ext_picard} (Picard), as compared with: SINDy; MATLAB's \texttt{nlgreyest} (w. fmincon) (fmincon); a black box gradient-free alternative approach to solving Eqn.~\eqref{eqn:optimization_problem_unconstrained} (GA); MATLAB's AMIGO2 Toolbox; Julia SciMLSensitivity Toolbox. For each algorithm, ARE and standard deviations for each parameter $\rho$ (Fig. 3a), $\sigma$ (Fig. 3b), $\beta$ (Fig. 3c) are computed over 10 datasets for each level of measurement noise variance -- $\varepsilon \in \{0.01, 0.05, 0.1, 0.2, 0.4\}$. True parameter values are $\rho = 2.8$, $ \sigma = 1$, $ \beta = 8/3$.  }}\label{fig:lorenz_results}
\end{figure*}

 \subsection{The Lorenz System: Chaos and Comparison.}
 \label{subsec:example_lorenz_multi_trajectories}

The goal of this last numerical test is to evaluate the accuracy of estimated parameters in the proposed algorithms as compared to existing state-of-the-art methods for parameter estimation. For this analysis we use data generated from the well-studied Lorenz system~\cite{brunton2022data, mauroy2016linear}, originally proposed as a model of convection rolls in the atmosphere~\cite{lorenz1963deterministic}. Solutions to the Lorenz system tend to converge to a 2D manifold known as a ``strange attractor'', but are chaotic within this manifold -- implying adjacent trajectories diverge exponentially fast. As a result, while parameter estimation for the Lorenz system is reasonably accurate when measurement noise is small, many methods for parameter estimation fail when noise is significant. The Lorenz system is defined as   \vspace{-3mm}
\begin{align}
\dot x_1(t) &= 10\sigma(x_2(t) - x_1(t)) \notag \\
\dot x_2(t) &= 10x_2(t)(\rho - x_3(t)) - x_2(t) \label{eqn:example_lorenz_model} \\
\dot x_3(t) &= 10x_1(t) x_2(t) - \beta x_3(t), \notag \\[-8mm]\notag
\end{align}
where $x_1, x_2, x_3$ are states and $\sigma, \rho, \beta$ are parameters. 

For this analysis, we use Lorenz parameters $\sigma = 1$, $\beta = \frac{8}{3}$ and $\rho = 2.8$. In order to allow for comparison with existing methods, each dataset uses a single initial condition, and directly measures the state at frequent, regular sampling times $t_i = \{0: 2i/(N_s - 1) :  2\}$ ($N_s = 101$). Measurement noise is modelled as $
y_i = (1 + n_i)\phi(t_i, x_0, \theta) + m_i,$
where $n_i, m_i$ are normally distributed with zero mean and variance $\varepsilon$. This approach is used to obtain  $50$ data sets $Y_{\varepsilon, x_0} = \{(y_i, t_i) \}_{i = 1}^{N_s}$, where $\varepsilon \in \{0, 0.01, 0.05, 0.1, 0.2, 0.4\}$ and 10  randomly chosen initial conditions $x_0$, normally distributed with zero mean and variance $1$.

Parameter estimation is then performed on each dataset using each of 6 algorithms: Picard (Alg.~\ref{alg4:ext_picard}); SINDy~\cite{brunton2022data}; fmincon~\cite{byrd1988approximate};  GA~\cite{sivanandam2008genetic}; { AMIGO2~\cite{balsa2016amigo2}; and SciML Sensitivity Toolbox~\cite{rackauckas2020universal}}. For Picard, Alg.~\ref{alg4:ext_picard} is used with parameters $\lambda = 5$, $n = 2$ and $T = 0.05$. For SINDy, we use a slight modification of the black-box implementation given in~\cite{brunton2022data} (regularization parameter $\lambda=0.25$ and measurements of $\dot x(t_i)$ are computed using first-order difference). For fmincon, we use \texttt{nlgreyest} from MATLAB's System Identification Toolbox~\cite{ljung2012version} with selected optimizer: \texttt{fmincon}. {  For AMIGO2 Toolbox, we use Enhanced Scatter Search to find the optimal parameters~\cite{egea2007scatter}. For SciML Sensitivity Toolbox, we use Adam optimizer followed by BFGS optimizer~\cite{fletcher2013practical}.} For GA, we designed a  gradient-free black-box optimization approach to solution of Eqn.~\eqref{eqn:optimization_problem_unconstrained}. Specifically, for given $x,\theta$, we use numerical simulation to compute the solution $\phi$ which then can be used to compute the least squares loss. With this loss oracle in hand, we use a genetic algorithm (MATLAB's \texttt{ga} function -- with tolerances $0.01$ and run time of $5$ mins on Intel i7-4960X CPU at 3.60 GHz) to solve the optimization problem. For all methods, parameters are restricted to $\rho \in [1.4 , 4.2]$, $\sigma \in [0.5 ,1.5]$, $\beta \in [\frac{4}{3} , 4]$, and  $x_i(0)=[0.5 x_i , 1.5 x_i]$, where $x_i$ are true initial conditions. 
  
\begin{figure}
 \centering
  
     \includegraphics[width = 0.65 \linewidth, trim = 10 40 40 40,clip]{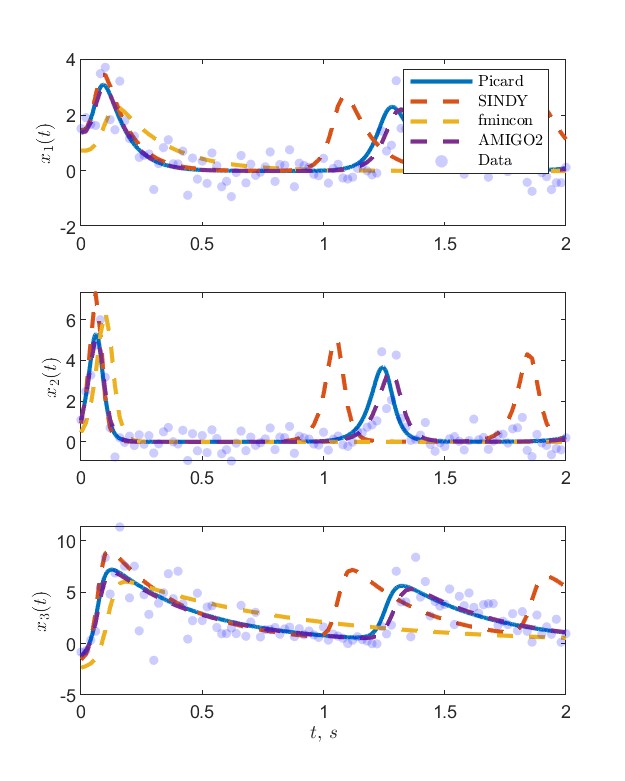}    \vspace{-3.5mm} 
   \caption{  {\footnotesize Numerical simulation of 3 states of the Lorenz system (Eqn.~\eqref{eqn:example_lorenz_model}) using identified parameters from Picard (Alg.~\ref{alg4:ext_picard}), SINDy, AMIGO2 and fmincon alternatives (See Fig.~\ref{fig:lorenz_results} and description) for noise variance $\varepsilon=0.4$. True parameter values are $\rho =2.8$, $\sigma = 1$, $\beta = 8/3$; initial states are $x(0)\hspace{-1mm} =\hspace{-1mm} [1.46,  0.90,  -1.57]^T$\hspace{-1mm}; and blue dots are measurements. }
}\vspace{-1.5mm}
    
 \label{fig:lorenz_results_example}
\end{figure}

Each of the four algorithms was evaluated on each of the 50 datasets and for each $\varepsilon$, we computed ARE and standard deviation as described previously and illustrated in Fig~\ref{fig:lorenz_results}. These results show that the Picard method outperforms all methods for low level of noise $\varepsilon \leq 0.1$ and for larger values of noise ($\varepsilon \geq 0.2$) Picard performs significantly better than SINDy and fmincon methods -- but comparable to the proposed black box optimization alternative to the gradient-based Picard approach. Simulation of the Lorenz system and comparison to data for estimated parameters in each method (with noise variance $\varepsilon=0.4$ and single dataset) is given in Fig.~\ref{fig:lorenz_results_example}.

\section{Conclusion}
In this paper, we reformulate the problem of parameter estimation for nonlinear ordinary differential equations as a constrained optimization problem with infinite-dimensional variables and constraints. We then propose a new class of gradient-contractive algorithms, based on the contractive properties of the Picard iteration, which eliminates the need for infinite-dimensional constraints and variables (Algs.~\ref{alg2:picard}-\ref{alg4:ext_picard}). In contrast with existing methods, this approach is gradient-based and allows for data which is sparse, irregular, and includes only partial measurements of the state and does not require numerical simulation or measurements of time derivatives. Furthermore, we have proposed sufficient conditions under which the proposed algorithm converges to a local optima. Finally, we have exhaustively tested several aspects of performance of the algorithm on a battery of nonlinear models and noisy datasets. 
\section*{Acknowledgments} We would like to thank Rolf Findeisen for his feedback and acknowledge the priority of his team in recognizing the role of Picard iteration in parameter estimation.

{
\bibliographystyle{plain}
\bibliography{Talitckii_Picard}
}
\appendix
%
\section{Auxiliary Lemmas}
\label{apd:proof_of_thm7} 

In this appendix, we prove the auxiliary lemmas.
\begin{lem}\label{lem:part_of_proof_of_thm7} 
Let $U \subset \mcl U:=\mcl C_0(\R,\R^{n_u})$ and $X \subset \R^{n_x}$ be convex and compact sets with $r=\sup_{\mbf u\in U}\norm{u}_{\mcl U}$.  Suppose that $L_i$ are differentiable in $x \in X$ with $\nabla_x L_i(x,u)$ Lipschitz on $(x,u) \in X \times B_r$.

 Let $\mcl P:X \times U \rightarrow U$ and there exists $N \N$ such that \vspace{-3.5mm}   
\begin{align}\label{eqn:general_picard_is_lipschitz_apdx}
&\|\mcl P^n(x_1, \mbf u_1) - \mcl P^n(x_2, \mbf u_2)\|_{\mcl U} \leq q^n_{\mcl P} \|\mbf u_1 - \mbf u_2\|_{\mcl U} \\
& \qquad \qquad \qquad  + K_{\mcl P} \|x_1 - x_2\|_X  \notag\\
&\|\nabla_x \mcl P^n(x_1, \mbf u_1) - \nabla_x \mcl P^n(x_2, \mbf u_2)\|_{\mcl U} \leq q^n_{\mcl P} \|\mbf u_1 - \mbf u_2\|_{\mcl U} \notag \\
& \qquad \qquad \qquad  + K_{\mcl P} \|x_1 - x_2\|_X \notag\\[-8mm] \notag
\end{align}
for all $n \geq N$, $x_1, x_2 \in X$ and $\mbf u_1, \mbf u_2 \in U$.

Then, there exist $K_L \ge 0$ and $K_{n, u}$, such that $\lim_{n \rightarrow \infty} K_{n, u} = 0$ and for all $n \in \N$, $x_1, x_2 \in X$ and $\mbf u_1, \mbf u_2 \in U$.  \vspace{-3.5mm}
\begin{align*}
&\|\nabla_x L(x_1, \mcl P^n(x_1, \mbf u_1)) -\nabla_x  L(x_2, \mcl P^n(x_2, \mbf u_2))\|_2 \\
&\qquad\qquad \leq K_{n,\mbf u} \|\mbf u_1 - \mbf u_2\|_{\mcl U} + K_{L} \|x_1 - x_2\|_X \\[-8mm]\notag
\end{align*}
where $L(x, \mbf u) = \sum_i L(x, \mbf u(t_i))$.
\end{lem}\vspace{-3mm}

\begin{pf} For notational clarity, we explicitly label the Lipschitz continuous partial gradients as $\nabla_1 L_i(x,v):=\nabla_x L_i(x,v)$ and $\nabla_2 L_i(x,v):=\nabla_v L_i(x,v)$. Then by the chain rule, for any  $n \in \N$, $x \in X$ and $\mbf u \in U$ we have  \vspace{-3.5mm}
\begin{align}
&\nabla_x L(x, \mcl P^n(x, \mbf u)) = \sum\nolimits_{i} \big[ \nabla_1 L_i(x, \mcl P^n(x, \mbf u)(z_i)) \label{eqn:nablax}
\\[-1mm]
 &\qquad\qquad\quad  + \nabla_2  L_i(x, \mcl P^n(x, \mbf u)(z_i))  \nabla_x \mcl P^n(x, \mbf u)(z_i)\big].\notag \\[-8mm]\notag
\end{align}

Note, since $X, U$ are compact sets and $L(x, \mbf u) = \sum_i L(x, \mbf u(t_i))$, $\nabla_{j} L_i(x, v)$  are bounded and  Lipschitz continuous functions of $x,v$ for $j \in \{1,2\}$. Thus, for some constant $K_2 > 0$, for all  $n \geq N$, $x_1, x_2 \in X$ and $\mbf u_1, \mbf u_2 \in U$ and $j \in \{1,2\}$ we have \vspace{-3.5mm} 
\begin{align*}
&\|\nabla_{j} L_i(x_1, \mcl P^n(x_1, \mbf u_1)) -  \nabla_{j} L_i(x_2, \mcl P^n(x_2, \mbf u_2))\|_2 \\ 
&\leq K_2 (\|x_1 - x_2\| + \|\mcl P^n(x_1, \mbf u_1)(z_i) - \mcl P^n(x_2, \mbf u_2)(z_i)\|_2) \\ 
&\qquad\leq q^n_{\mcl P}K_2 \|\mbf u_1 - \mbf u_2\|_{\mcl U} + K_2(1+K_{\mcl P})\|x_1 - x_2\|_2. \\[-8mm]
\end{align*}

Therefore, Condition~\eqref{eqn:general_picard_is_lipschitz_apdx} and Eqn.~\eqref{eqn:nablax} imply that  $\nabla_x L(x, \mcl P^n(x, \mbf u))$ is Lipschitz continuous with respect to $x$ and $\mbf u$ (as the product and sum of Lipschitz and bounded functions) and where the Lipschitz factor with respect to $\mbf u$ decreases with increasing $n$. Specifically, there exist $K_L \geq 0$ and $\{K_{n,\mbf u}\}_{n = 0}^{\infty} \subset \R_+ $  such that $\lim_{n\rightarrow \infty} K_{n,\mbf u} = 0$ and \vspace{-3.5mm}
\begin{align*}
&\|\nabla_x L(x_1, \mcl P^n(x_1, \mbf u_1)) -\nabla_x  L(x_2, \mcl P^n(x_2, \mbf u_2))\|_2 \\
&\qquad\qquad \leq K_{n,\mbf u} \|\mbf u_1 - \mbf u_2\|_{\mcl U} + K_{L} \|x_1 - x_2\|_X \\[-8mm]\notag
\end{align*}
for all $n \geq N$, $x_1, x_2 \in X$ and $\mbf u_1, \mbf u_2 \in U$. $\blacksquare$

\end{pf}\vspace{-2mm}

\section{Gradients of Picard Iterations}
\label{apd:Sensitivity} 

In this appendix, we show that gradients of Picard Iterations are Lipschitz continuous. Specifically, we first show that $\nabla_{x, \theta} \mcl P_{t_0, x, \theta} \mbf u$ is Lipschitz continuous with respect to $\mbf u$ uniformly in $(x, \theta) \in X\times \Theta$.
\begin{lem}  \label{lem:picard_sensitivity_is_continuous_in_phi}
Suppose the conditions of Lem.~\ref{lem:picard_sensitivity_is_cauchy} are satisfied. 
Then, there exists $K_{s} \geq 0$ such that  for any $n \geq m$, $(x, \theta) \in X\times \Theta$  and $\mbf u, \mbf v \in C_a[0, T]$  we have \vspace{-3.5mm} 
\begin{align*}
& \|\nabla_{x, \theta} {\mcl P}_{t_0, x, \theta}^n \mbf u - \nabla_{x, \theta} {\mcl P}_{t_0, x, \theta}^m \mbf v\|_\infty  \\
& \quad\quad\leq (TK_{x})^{m} \bigg[ mK_{s} \|{\mcl P}_{t_0, x, \theta}^{n-m}\mbf u- \mbf v\|_\infty +  \| \nabla_{x, \theta}  \mcl P^{n-m}_{t_0, x, \theta} \mbf u \|_\infty \bigg]. \\[-8mm]
\end{align*} 
\end{lem}\vspace{-5mm}
\begin{pf}  
Define composite functions $\mbf u_n(t, x, \theta) := ({\mcl P}_{t_0, x, \theta}^n \mbf u)(t)$ and $\mbf v_m(t, x, \theta) = ({\mcl P}_{t_0, x, \theta}^m \mbf v)(t)$. Then, for all $x \in X$, $\theta\in\Theta$ and $t\in[0, T)$ \vspace{-3.5mm} 
\begin{align*}
&\|\nabla_x \mbf u_n(t,   x, \theta) - \nabla_x \mbf v_m(t,  x, \theta)\|_2 \\
& = \bigg\|\int_{0}^t \hspace{-1mm}\nabla_x f(t_0 + s, \mbf u_{n-1}, \theta) \nabla_x \mbf u_{n-1}(s, x, \theta) \\
&\qquad\qquad\qquad - \nabla_x f(t_0 + s, \mbf v_{m-1}, \theta) \nabla_x \mbf v_{m-1}(s, x, \theta)ds\bigg\|_2. 
\end{align*} 
Note, the Lipschitz continuity of $f$, $\mbf u_n$, $\mbf v_n$ (from Lem.~\ref{lem:picard_is_lipschitz_in_parameters}) implies  $f$, $\nabla_{x, \theta} \mbf u_n$, $\nabla_{x, \theta} \mbf v_n$ are bounded.

%
Next, for some $K_1 \geq 0$ we have for all $x \in X$, $\theta\in\Theta$ and $t\in[0, T)$ \vspace{-3.5mm}
\begin{align*}
 & \|\nabla_x \mbf u_n(t, x, \theta) - \nabla_x \mbf v_m(t, x, \theta)\|_2 \\
 & \leq  K_{x} \left\|\int_{0}^t \nabla_x  \mbf u_{n-1}(s, x, \theta) - \nabla_x  \mbf v_{m-1}(s, x, \theta)ds\right\|_2 \\
 & + K_1 \left\|\int_{0}^t \nabla_x f(t_0 + s, \mbf u_{n-1}, \theta) -\nabla_x f(t_0 + s, \mbf v_{m-1}, \theta) ds\right\|_2\hspace{-2mm}.
\end{align*} 
Thus, for all $x \in X$, $\theta\in\Theta$ and $t\in[0, T)$ \vspace{-3.5mm}
\begin{align*}
&\sup_{t\in[0, T)}\|\nabla_x \mbf u_n(t,  x, \theta) - \nabla_x \mbf v_m(t,  x, \theta)\|_2 \\
& \leq  TK_{x} \bigg(\sup_{t\in[0,T)}\|\nabla_x \mbf u_{n-1}(t, x, \theta) - \nabla_x  \mbf v_{m-1}(t, x, \theta)\|_2  \\
& \qquad\qquad + K_1/K_x \sup_{t\in[0, T)}\|\mbf u_{n-1}(t,  x, \theta) -   \mbf v_{m-1}(t,   x, \theta)\|_2\bigg).
\end{align*}
Therefore, by induction we have for all $x\in X$ and $\theta\in\Theta$ we have \vspace{-3.5mm}
\begin{align*}
&\|\nabla_x {\mcl P}_{t_0, x, \theta}^n \mbf u)- \nabla_x {\mcl P}_{t_0, x, \theta}^m \mbf v\|_\infty \\
& \quad\leq (TK_{x})^{m} \bigg[ m K_{s} \|\mcl P^{n-m}_{t_0, x, \theta}\mbf u - \mbf v\|_\infty  +  \| \nabla_x  \mcl P^{n-m}_{t_0, x, \theta} \mbf u \|_\infty \bigg]. 
\end{align*} 
$\nabla_{\theta} \mcl P^n_{t_0, x, \theta} \mbf u$ can be considered analogously.$\blacksquare$   \vspace{-2.5mm}
 
\end{pf}\vspace{-2mm}

 Next, we show that $\nabla_{x, \theta} \mcl P_{t_0, x, \theta} \mbf u$ is Lipschitz continuous with respect to $x$ and $\theta$ uniformly in $\mbf u \in C_a[0, T]$.
\begin{lem}\label{lem:picard_sensitivity_is_continuous_in_p}
Suppose the conditions of Lem.~\ref{lem:picard_sensitivity_is_cauchy} are satisfied. Then,  there exists $K > 0$ such that  for any $n\in \N$, $(x_1, \theta_1), (x_2, \theta_2) \in X\times \Theta$  and $\mbf u \in C_a[0, T]$  we have \vspace{-3.5mm}
\begin{align*}
& \|\nabla_{x, \theta} {\mcl P}_{t_0, x_1, \theta_1}^n \mbf u - \nabla_{x, \theta} {\mcl P}_{t_0, x_2, \theta_2}^n \mbf u\|_\infty  \\
& \qquad\qquad\leq K \|x_1 - x_2\|_2 + K\|\theta_1 - \theta_2\|_2.
\end{align*} 
\end{lem}
\begin{pf}   
Define function $\mbf u_n(t, x, \theta):= ({\mcl P}^n_{t_0, x, \theta}\mbf u)(t)$, then for all $(x_1, \theta_1), (x_2, \theta_2) \in X\times \Theta$ and $\mbf u \in C_a[0, T]$\vspace{-3.5mm}
\begin{align*}
&\|\nabla_x \mbf u_n(t,x_1, \theta_1) - \nabla_x \mbf u_n(t, x_2, \theta_2)\|_2  \\
& =  \bigg\|\int_{0}^t \nabla_x f(t_0 + s, \mbf u_{n-1}, \theta_1) \nabla_x \mbf u_{n-1}(s, x_1, \theta_1) \\
&\hspace{15mm}- \nabla_x f(t_0 + s, \mbf u_{n-1}, \theta_2) \nabla_x \mbf u_{n-1}(s, x_2, \theta_2)ds\bigg\|_2  
\end{align*}
Next, since $\nabla_x f$,  $\nabla_x \mbf u_{n-1}(s, x_2, \theta_2)$ are bounded (since $f$ and $\mbf u_{n-1}$ are Lipschitz from Lem.~\ref{lem:picard_is_lipschitz_in_parameters}) there exists $K_1 \geq 0$ such that\vspace{-3.5mm}
\begin{align*}
&\|\nabla_x \mbf u_n(t,x_1, \theta_1) - \nabla_x \mbf u_n(t, x_2, \theta_2)\|_2  \\
 & \leq  K_{x} \left\|\int_{0}^t \nabla_x  \mbf u_{n-1}(s, x_1, \theta_1) - \nabla_x  \mbf u_{n-1}(s, x_2, \theta_2)ds\right\|_2 \\
 & \qquad\qquad\qquad + K_1 \bigg\|\int_{ 0}^t \nabla_x f(t_0 + s, \mbf u_{n-1}(t,x_1, \theta_1), \theta_1) \\
 & \qquad\qquad\qquad\qquad\quad-\nabla_x f(t_0 + s, \mbf u_{n-1}(t,x_2, \theta_2), \theta_2) ds\bigg\|_2 
\end{align*}
 for all  $(x_1, \theta_1), (x_2, \theta_2) \in X\times \Theta$ and $\mbf u \in C_a[0, T]$.
 
Thus, since $\nabla_x f$ and $\mbf u_{n-1}$ is Lipschitz continuous in $x$ and $\theta$ there exists $K_2\geq 0$ such that \vspace{-3.5mm}
\begin{align*}
&\sup_{t\in[0, T)}\|\nabla_x \mbf u_n(t,  x_1, \theta_1) - \nabla_x \mbf u_n(t,  x_2, \theta_2)\|_2 \\
&\leq TK_{x}  \sup_{t\in[0, T)}  \|\nabla_x   \mbf u_{n-1}(t,  x_1, \theta_1) - \nabla_x  \mbf u_{n-1}(t,  x_2, \theta_2)\|_2  \\
 & \qquad\qquad\qquad\qquad +K_2 \|x_1 - x_2\| + K_2 \|\theta_2 - \theta_1\|_2. 
 \end{align*} 
 for all  $(x_1, \theta_1), (x_2, \theta_2) \in X\times \Theta$ and $\mbf u \in C_a[0, T]$.
 
 Next, since $\mbf u \in C_a[0, T]$ is not a function of $x$ and $\theta$, by induction for all  $(x_1, \theta_1), (x_2, \theta_2) \in X\times \Theta$ and $\mbf u \in C_a[0, T]$ we have \vspace{-3.5mm}
\begin{align*}
&\sup_{t\in[0, T)}\|\nabla_x \mbf u_n(t,  x_1, \theta_1) - \nabla_x \mbf u_n(t,  x_2, \theta_2)\|_2 \\
& \leq \frac{K_2}{1 - TK_x}\|x_1 - x_2\| + \frac{K_2}{1 - TK_x} \|\theta_2 - \theta_1\|_2. \\
\end{align*} 
Analogously, we have that $\nabla_p \mbf u_n(t, x, p)$ is Lipschitz continuous. Thus, for $K=\frac{K_2}{1 - TK_x}$ we have \vspace{-3.5mm}
\begin{align*}
& \|\nabla_{x, \theta} {\mcl P}_{t_0, x_1, \theta_1}^n \mbf u - \nabla_{x, \theta} {\mcl P}_{t_0, x_2, \theta_2}^n \mbf u\|_\infty  \\
& \qquad\qquad\leq K \|x_1 - x_2\|_2 + K\|\theta_1 - \theta_2\|_2
\end{align*}
  for any $n\in \N$, $(x_1, \theta_1), (x_2, \theta_2) \in X\times \Theta$  and $\mbf u \in C_a[0, T]$.   $\blacksquare$  
\end{pf}\vspace{-2mm}

\end{document}